\documentclass[a4paper,10pt]{article} \usepackage{amsmath,amssymb,amsthm}
\usepackage[english]{babel}

\newtheorem{theo}{Theorem}[section] \newtheorem{defi}[theo]{Definition}
\newtheorem{lemm}[theo]{Lemma} \newtheorem{prop}[theo]{Proposition}
\newtheorem{coro}[theo]{Corollary}

\newcommand{\Na}{\mathbb N}                   

\newcommand{\Ra}{\mathbb R}                   

\newcommand{\Ca}{\mathbb C}                   

\newcommand{\scal}[1]{\langle #1 \rangle}

\newcommand{\finpreuve}{\hfill $\Box$}

\newcommand{\name}{$\underline{\qquad \qquad}$}

\newcommand{\refe}[1]{\ref{#1}} \newcommand{\reff}[1]{(\ref{#1})}

\begin{document}

\author{  Jean-Marc
Bouclet\footnote{Jean-Marc.Bouclet@math.univ-lille1.fr}\\
 \\ Universit\'e  de Lille  1 \\ Laboratoire Paul Painlev\'e \\ UMR  CNRS 8524,  \\ 59655 Villeneuve  d'Ascq }
\title{{\bf \sc  Resolvent estimates for the Laplacian on asymptotically hyperbolic manifolds}}

\maketitle

\begin{abstract} Combining results of Cardoso-Vodev \cite{CaVo1} and Froese-Hislop \cite{FrHi1},
we use Mourre's theory to prove high energy estimates for the
 boundary values of the weighted resolvent of the Laplacian on an
asymptotically hyperbolic manifold. We derive estimates involving
a class of pseudo-differential weights which are more natural in
the asymptotically hyperbolic geometry than the weights $
\scal{r}^{-1/2-\epsilon} $ used in \cite{CaVo1}.
\end{abstract}

\section{Introduction, results and notations}
\setcounter{equation}{0}

The purpose of this paper is to prove resolvent estimates for the
Laplace operator $ \Delta_g $ on a non compact Riemannian manifold
$ ({\mathcal M},g) $ of asymptotically hyperbolic type. The latter
means that $ {\mathcal M} $ is a connected manifold of dimension
$n$ with or without boundary such that, for some relatively
compact open subset $ {\mathcal K} $, some closed manifold $ Y $
(i.e. compact, without boundary) and some $ r_0
> 0 $, $ ({\mathcal M} \setminus {\mathcal K},g) $ is isometric to
$ [r_0,+\infty) \times Y $ equipped with a metric of the form
\begin{eqnarray}
dr^2 + e^{2r} h (r) .  \label{metrique}
\end{eqnarray}
For   each $r$, $ h (r) $ is a
Riemannian metric on $ Y $ which is a perturbation of a fixed
metric $ h $, meaning that, for all $ k $ and all semi-norm $
|||.||| $ of the space of smooth sections of $  T^* Y \otimes T^*
Y $,
\begin{eqnarray}
 \sup_{r \geq
r_0 } \left| \left| \left| \scal{r}^2 \partial^k_r ( h (r) -
h ) \right| \right| \right| < \infty, \label{decayperturb}
\end{eqnarray}
with $ \scal{r} = (1 + r^2)^{1/2} $. Here, and in the sequel, $ r
$ denotes a positive smooth function on $ {\mathcal M} $ going to
$ + \infty $ at infinity and  which is a coordinate near $
{\mathcal M} \setminus {\mathcal K} $, i.e. such that $ dr $
doesn't vanish near $ {\mathcal M} \setminus {\mathcal K} $. Such
manifolds include the hyperbolic space $ \mathbb{H}_n $ and some
of its quotients by discrete isometry groups. More generally, we
have typically in mind  the context of the $0$-geometry of Melrose
\cite{Melr0}.

Let $ G $ be the Dirichlet or Neumann realization of $ \Delta_g $
(or the standard one if $
\partial M $ is empty) on $ L^2 ({\mathcal M},d\mbox{Vol}_g) $. Then, according
to \cite{CaVo1}, it is known that  the limits $ \scal{r}^{-s} (G -
\lambda \pm i 0)^{-1} \scal{r}^{-s} := \lim_{\varepsilon
\rightarrow 0^+} \scal{r}^{-s} (G - \lambda \pm i \varepsilon
)^{-1} \scal{r}^{-s} $ exist, for all $ s
> 1/2 $, and satisfy
\begin{eqnarray}
\left| \left| \scal{r}^{-s} (G - \lambda \pm i 0)^{-1}
\scal{r}^{-s} \right| \right|_{L^2({\mathcal M}, d
\mathrm{Vol}_g)} \leq C e^{C_G \lambda^{1/2}}, \qquad \lambda \gg
1 . \label{CardosoVodev}
\end{eqnarray}
In \cite{Vode2}, it is shown that the right hand side can be
replaced by $ C \lambda^{-1/2} $, under a non trapping condition.

In the present paper, we will mainly prove that, up to logarithmic
terms in $ \lambda $, such estimates still hold if one replaces $
\scal{r}^{-s} $ by a class of operators which are, in some sense,
weaker than $ \scal{r}^{-s} $ and more adapted to the framework of
the asymptotically hyperbolic scattering.

\medskip

Let us fix the notations used in this article.

Throughout the paper, $ C_c^{\infty} ({\mathcal M}) $ denotes the
space of smooth functions with compact support. If $ {\mathcal M}
$ has a boundary, $ C_0^{\infty} ({\mathcal M}) $ is the subspace
of $ C_c^{\infty}({\mathcal M}) $ of functions vanishing near $
\partial {\mathcal M} $ and if $ B $ denotes the boundary conditions
associated to $ G $ (if any), $ C_B^{\infty} ({\mathcal M}) $ is
the subspace of $ \varphi \in C_c^{\infty}({\mathcal M}) $ such
that $ B \varphi = 0 $ (e.g. $ B \varphi = \varphi_{|\partial M}
$ for the Dirichlet condition).

 We set $ I = (r_0 , + \infty) $ and call $ \iota $ the
isometry from $ {\mathcal M} \setminus {\mathcal K} $ to $ \bar{I}
\times Y $. If $  \Psi : U_Y \subset Y \ni \omega \mapsto
(y_1,\cdots,y_{n-1}) \in U \subset \Ra^{n-1}  $ is a coordinate
chart and $ {\mathcal M} \setminus {\mathcal K} \ni m \mapsto
\omega (m) \in Y $ is the natural projection induced by $ \iota $,
we define the chart  $ \tilde{\Psi} : \iota^{-1} ( I \times U_Y )
\subset {\mathcal M} \rightarrow I \times U $   by
\begin{eqnarray}
 \tilde{\Psi} (m) = \left( r (m), \Psi (\omega(m))  \right) .
\label{chart}
\end{eqnarray}
There clearly exists a finite atlas on $ {\mathcal M} $ composed
of such charts and compactly  supported ones. For any diffeomorphism $f: M
\rightarrow N  $,
between open subsets of two manifolds, we use the standard notations $ f^* 
 $ and $ f_*  $
for the maps defined by $f^* u = u \circ f^{-1}  $ and $ f_* u = u \circ f $,
respectively on $ C^{\infty}(M) $ and $ C^{\infty}(N) $ (and more generally on
differential forms or sections of density bundles). 

By $ \reff{metrique} $ and $ \reff{decayperturb} $,  we have $  \iota^* (  d
\mbox{Vol}_g ) = \tilde{\Theta}   e^{(n-1)r}
drd\mbox{Vol}_h  $ on $ {\mathcal M} \setminus {\mathcal K} $,
with $ \tilde{\Theta} = d \mbox{Vol}_{h(r)} / d \mbox{Vol}_h $
satisfying $ \sup_{I} ||| \scal{r}^2
\partial_r^k ( \tilde{\Theta} (r,.) - 1) ||| < \infty  $ for all $ k $ and
all seminorm $ |||. ||| $ of $ C^{\infty}(Y) $. We  choose a
positive function $ \Theta \in C^{\infty}({\mathcal M}) $ such
that $ \iota^*  \Theta  =e^{(n-1)r} \tilde{\Theta} $ on $
{\mathcal M} \setminus {\mathcal K} $ and we define a new measure
$ d\mbox{Vol}_{\mathcal M} = \Theta^{-1} d\mbox{Vol}_g $. This is
convenient since we now have $ \iota^* ( d\mbox{Vol}_{\mathcal M} ) 
= drd\mbox{Vol}_h   $ on $ I \times Y  $
hence, if we set
   $ L^2 (\mathcal M) =  L^2 (\mathcal M ,
d\mbox{Vol}_{\mathcal M}  ) $, we get natural unitary isomorphisms
\begin{eqnarray}
L^2 ({\mathcal K}) \oplus L^2 ({\mathcal M}
\setminus {\mathcal K})  \approx L^2 ({\mathcal K}) \oplus
L^2 (I , dr ) \otimes L^2 (Y, d\mbox{Vol}_h) \approx  L^2 ({\mathcal K} ) \oplus
\bigoplus_{k=0}^{\infty} L^2 (I,dr) , \label{unitary}
\end{eqnarray}
using, for the last one, an
orthonormal basis $ (\psi_k)_{k \geq 0} $ of eigenfunctions of $ \Delta_h
$. More explicitely,  the isomorphism between $ L^2 (I , dr ) \otimes L^2 (Y,
d\mbox{Vol}_h)  $ and $
\bigoplus_{k=0}^{\infty} L^2 (I,dr)  $ is given by $ \varphi
\mapsto (\varphi_k)_{k \geq 0} $ with
\begin{eqnarray}
\varphi_k (r) = \int_Y \varphi (r,\omega)
\overline{\psi_k(\omega)} \ d\mbox{Vol}_h (\omega) .
\label{coefFourier}
\end{eqnarray}

In what follows, we will consider the self-adjoint operator
$$ H = \Theta^{1/2} G \Theta^{-1/2} $$ on $ L^2 (\mathcal M) $, with
domain $ \Theta^{1/2} D (G)  $. If $ \partial {\mathcal M} $ is
non empty, we furthermore assume that $ \Theta \equiv 1 $ near $
\partial {\mathcal M} $ in order to preserve the boundary condition. This is an elliptic differential
operator, unitarily equivalent to $ G $, which takes the form, on
$ {\mathcal M} \setminus {\mathcal K}  $,
\begin{eqnarray}
 H =  D_r^2 + e^{-2r} \Delta_h + V  + (n-1)^2/4  \label{notationcarte},
\end{eqnarray}
with $ \Delta_h $ the Laplace operator on $ Y $ associated to the $r$-independent
metric $h$ and  $ V $  a second order differential operator of the
following form in local coordinates
\begin{eqnarray}
 \tilde{\Psi}^* V  \tilde{\Psi}_* = \sum_{   |\beta| \leq 2 }
\scal{r}^{-2} v_{\beta} (r,y) (e^{-r} D_y)^{\beta}  ,
\label{formeV}
\end{eqnarray}
 with $
\partial_r^k \partial_y^{\alpha} v_{\beta} $ bounded on $ I
\times U_0 $ for all $ U_0 \Subset  U $ and all $ k , \alpha $.
Here $ U $ is associated to the chart $ \Psi $ (see  above $ \reff{chart}$).
 Without loss of generality, by possibly
increasing $ r_0 $, we may assume that
$$ H = H_0 + V $$
 with $ V $
 of the same form  as above, with coefficients supported in $
{\mathcal M} \setminus {\mathcal K} $, which is $ H $  bounded
with relative bound $ < 1 $ (see Lemma $1.4$ of \cite{FrHi1} or
Lemma $ \refe{FrHiSobolev} $ below), and $ H_0 $ another
self-adjoint operator (with the same domain as $H$) such that
\begin{eqnarray}
 H_0 =  D_r^2 + e^{-2r} \Delta_h   + (n-1)^2/4, \label{Hzero}
\end{eqnarray}
on  $ \iota^{-1} \left( (r_0 + 1,\infty ) \times Y \right) $.

We next choose a positive function $ w \in C^{\infty}(\Ra) $ such that
\begin{eqnarray}
 w (x) =  \begin{cases}
    1, & x \leq 0, \\
    x, & x \geq 1
  \end{cases} . \label{echelleclassique}
\end{eqnarray}
If $ \mbox{spec} (\Delta_h )  =  (\mu_k)_{k \geq 0} $ and $ s \geq
0 $, we define a bounded operator $ \widetilde{W}_{-s} $ on $ L^2 (I)
\otimes L^2 (Y,d\mbox{Vol}_h) $  by
\begin{eqnarray}
 (\widetilde{W}_{-s} \varphi) (r,\omega) = \sum_{k \geq 0} w^{-s}
(r -  \log \sqrt{ \scal{\mu_k} } ) \varphi_k (r) \psi_k (\omega) .
\label{optimum}
\end{eqnarray}
 Using $ \reff{unitary}  $, we pull $ \widetilde{W}_{-s} $ back as
an operator $ W_{-s} $ on $ L^{2}({\mathcal M}) $, assigning   $
W_{-s} $ to be the identity on ${L^2 (\mathcal K)}  $. We can now
state our main result.
\begin{theo} \label{theo1} Assume that, for some function $ \varrho (\lambda) \geq
c  \lambda^{-1/2} $ and some real number $ 0 < s_0  \leq 1 $,
\begin{eqnarray}
 || \scal{ r }^{-s_0} (H- \lambda \pm i 0)^{-1}  \scal{ r }^{-s_0}  || \leq C \varrho
(\lambda), \qquad \lambda \gg 1 . \label{abstraitdecroissance}
\end{eqnarray}
Then, for all $ s > 1/2  $, there exists $ C_s $ such that
\begin{eqnarray}
 || W_{-s} (H- \lambda \pm i 0)^{-1}  W_{-s} || \leq C_s (\log
\lambda)^{ 2 s_0 + 2  s } \varrho (\lambda)   , \qquad \lambda \gg
1 . \label{refpoids}
\end{eqnarray}
\end{theo}
Using the results of \cite{CaVo1,Vode2}, i.e. the estimates $ \reff{CardosoVodev} $, we obtain
\begin{coro} Let $ W^{\Theta}_{-s} = \Theta^{-1/2} W_{-s} \Theta^{1/2} $ with $s > 1/2 $. On any asymptotically hyperbolic manifold, we have
$$  \left| \left| W_{-s}^{\Theta} (G- \lambda \pm i 0)^{-1}  W_{-s}^{\Theta} \right| \right|_{L^2({\mathcal M}, d
\mathrm{Vol}_g)} \leq C_s (\log \lambda)^{4 s} e^{ C_G \lambda^{1/2}} , \qquad \lambda \gg 1 , $$ with the same $ C_G $
as in $ \reff{CardosoVodev} $. If the manifold is non trapping (in
the sense of \cite{Vode2}), we have
$$  \left| \left| W_{-s}^{\Theta} (G- \lambda \pm i 0)^{-1}  W_{-s}^{\Theta}
\right| \right|_{L^2({\mathcal M}, d \mathrm{Vol}_g)} \leq C (\log
\lambda)^{4 s}  \lambda^{-1/2}   , \qquad \lambda \gg 1 . $$
\end{coro}

These results improve the estimate $ \reff{CardosoVodev} $ to the
extent that $ W_{-s} $ and  $ W_{-s}^{\Theta}  $ are "weaker" than $
\scal{r}^{-s} $ in the sense that $ W_{-s} \scal{r}^s $ is not
bounded. The latter is easily verified using $ \reff{optimum} $ by
choosing a sequence $ (\varphi_k)_{k \geq 0} \in L^2 (I) $ such
that $ \sum_k ||\varphi_k||^2 = 1 $ with $ \varphi_k $ supported
close to $ \log \sqrt{ \scal{\mu_k} } $.

A result similar to Theorem $ \refe{theo1} $ has already been
proved by Bruneau-Petkov in \cite{BrPe1} for Euclidean scattering
(on $\Ra^n$). They essentially show that, if $ P $ is a long range
perturbation of  $- \Delta_{ \Ra^n } $ such that $ ||\chi
(P-\lambda \pm i 0)^{-1} \chi || = {\mathcal O}(e^{C \lambda}) $
for all $ \chi \in C_0^{\infty}(\Ra^n) $, then $ || \scal{x}^{-s}
(P-\lambda \pm i 0)^{-1} \scal{x}^{-s} || = {\mathcal O}(e^{C_1
\lambda})  $, with $ s > 1/2 $. In other words,  one can replace
compactly supported weights by polynomially decaying ones.

Weighted resolvent estimates can be used for various applications among which are spectral
 asymptotics, analysis of scattering matrices, of scattering amplitudes 
or non linear problems. In particular,
they are known to be useful to obtain Weyl formulas for
scattering phases in Euclidean scattering
\cite{Robe1,Robe2,BrPe1,Bouc4} and the present paper was motivated
by similar considerations in the hyperbolic context \cite{Bouc5,Bouc8}.
Actually, high energy estimates are
important tools to get semiclassical approximations of the
Schr\"odinger group by the techniques of Isozaki-Kitada \cite{IsKi1,IsKi2}. 
This is well known on $ \Ra^n $ \cite{Robe1,Robe2,Bouc4} and is being
developed for asymptotically hyperbolic manifolds \cite{Bouc5,Bouc8}.
These applications will be published elsewhere (they would otherwise lead to a
 paper of unreasonnable length).

\medskip

We now introduce  a class of pseudo-differential operators
associated with the scale of weights defined by the operators $
W_{-s} $. For $s \in \Ra$, we set 
$$ w_s (r,\eta) = w^{s} (r - \log \scal{\eta}) $$ 
and define the space $ {\mathcal S} (w_s) \subset
C^{\infty}(\Ra_r \times
 \Ra^{n-1}_y \times \Ra_{\rho} \times \Ra^{n-1}_{\eta})  $ as the
 set of
symbols satisfying
$$ \left| \partial_r^j \partial_y^{\alpha} \partial_{\rho}^k \partial_{\eta}^{\beta}
a(r,y,\rho,\eta) \right| \leq C_{j\alpha k \beta} w_s \left( r ,
\eta \right) , \qquad r,\rho \in \Ra, \ \ y,\eta \in \Ra^{n-1}.
$$
Note that, $ w_s  $ is a temperate weight in the sense of
\cite{Horm3} (see Lemma $ \refe{poids} $ of the present paper).
Note also that $ {\mathcal S}(w_{s_1}) \subset {\mathcal S}(w_{s_2}) $
if $ s_1 \leq s_2 $.

To construct operators on the manifold $ {\mathcal M} $, we consider a 
chart $ \Psi : U \rightarrow U_Y $ (we keep the notations 
above $ \reff{chart} $) and we choose open sets 
$ U_0 \Subset U_1 \Subset U_2 \Subset U  $. We pick cutoff functions 
$ \kappa , \tilde{\kappa} \in
C^{\infty}(\Ra_r \times \Ra_y^{n-1})
 $ which are respectively supported in $ I \times U_1 $ and $ I \times U_2
 $, with bounded derivatives and
 such that $ \tilde{\kappa} \equiv 1 $ near $ \mbox{supp} \ \kappa $,
$ \kappa \equiv 1 $ on $ (r_1 , + \infty) \times U_0 $ for some $
r_1 > r_0 $. For bounded symbols $a$, we can then 
define
$$   \tilde{\Psi}_* \kappa O \! p (a) \tilde{\kappa} \tilde{\Psi}^*
 = \tilde{\Psi}_* \kappa a(r,y,D_r,D_y) \tilde{\kappa} \tilde{\Psi}^*,
  $$
on $ L^2 (\mathcal M) $.
\begin{theo} \label{theo2} Assume that $ a \in {\mathcal S}(
w_{-s}) $ for some $ s \geq 0 $. Then, there exist bounded
operators $ B_{1,s} $ and $ B_{2,s} $ on $ L^2 ({\mathcal M}) $
such that
$$    \tilde{\Psi}_* \kappa O \! p (a) \tilde{\kappa} \tilde{\Psi}^* = B_{1,s} W_{-s} = W_{-s} B_{2,s} . $$
\end{theo}
The interest of this theorem is that Theorem $ \refe{theo1} $
still holds if one replaces $ W_{-s} $ by pseudo-differential
operators with symbols in $ {\mathcal S}(w_{-s})  $, $ s > 1/2 $. This is 
important since the classes $ {\mathcal S}(w_{-s})  $, with $ s > 0 $, are 
naturally associated with the functional calculus of asymptotically hyperbolic 
Laplacians as we shall see below.

Let us explain why polynomial weights $ \scal{r}^{-s} $ are more
natural for Euclidean scattering than for the asymptotically
hyperbolic one. In polar coordinates on $ \Ra^n $, the principal
symbol of the flat Laplacian is $ \rho^2 + r^{-2} q_0 $ (with $
q_0 = q_0 (y,\eta) $ the principal symbol of the Laplacian on the
sphere) and since $ d r^{-2}/dr = - 2 r^{-2} \times r^{-1} $, 
it is easy to check that, for all $ k \in \Na $, $
\gamma \in \Na^{n-1} $ and $ z \notin [0,+\infty) $, one has
\begin{eqnarray}
\left| \partial_r^k \partial_{\eta}^{\gamma}(\rho^2 +
r^{-2}q_0-z)^{-1} \right| \leq C_{z,k,\gamma}  |\rho^2 + r^{-2}q_0
-z|^{-1} r^{-k-|\gamma|}. \label{decayeuclidien}
\end{eqnarray}
Here we consider the function $ (\rho^2 + r^{-2}q_0 - z)^{-1} $
for it is the principal symbol of $ (-\Delta_{\Ra^n}-z)^{-1} $ (in
polar coordinates) and hence the prototype of the symbols involved
in
 the functional calculus of perturbations of $ -
\Delta_{\Ra^n} $. Besides, we note that when one considers a perturbation 
of $- \Delta_{\Ra^n} $ by a long range potential $V_L $, one usually assumes  
that, for some $ \varepsilon > 0 $,
$$ |\partial^{\alpha}_x V_L (x)| \leq C_{\alpha} \scal{x}^{-\varepsilon-
|\alpha|} . $$
Hence, powers of $ r^{-1} $ are naturally involved in the symbol classes for
Euclidean scattering. This is compatible with the fact that the  weights needed to
get resolvent estimates in this context are also powers of $ r^{-1}  $.

 In hyperbolic scattering, the situation is different. 
The principal symbol of $ H_0 $ (see $\reff{Hzero} $) takes the
form $ \rho^2 + e^{-2r} q_h $ (with $ q_h = q_h (y,\eta) $ the
principal symbol of $ \Delta_h $ on $Y$) and since $ d e^{-2r}/d
r = -2 e^{-2r} $ we cannot hope to gain any extra decay of the symbols with
respect to $ r $, unlike in the Euclidean case.  However, remarking that
$$ \left| \frac{e^{-2r}q_h}{\rho^2 + e^{-2r}q_h + 1} \right| \leq C_s
w_{-s}(r,\eta), \qquad \forall \ s \geq 0,
 $$
it is easy to check that, if $ k + |\gamma| \geq 1  $,
\begin{eqnarray}
\left| \partial_r^k \partial_{\eta}^{\gamma}(\rho^2 + e^{-2r
}q_h-z)^{-1} \right| \leq C_{z,k,\gamma,s}  |\rho^2 + e^{-2r}q_h
-z|^{-1} w_{-s}(r,\eta), \qquad \forall \ s \geq 0 .
\label{decayhyperbolique}
\end{eqnarray}
 Here again, we have
chosen 
$ (\rho^2 +
e^{-2r}q_h - z)^{-1}  $ since it is the principal symbol of the pseudo-differential
approximation of $ (H_0-z)^{-1}  $ (see \cite{Bouc5,Bouc8}). The estimate 
$ \reff{decayhyperbolique} $  reflects the fact that the weights $ w_{-s} $ are more
natural than $ \scal{r}^{-s} $ in hyperbolic scattering: we do not gain any
power of $ r^{-1} $ by differentiating but we gain powers of $ w_{-1} $ and 
these weights are naturally associated with the resolvent estimates as shown
by Theorems $ \refe{theo1}  $ and $ \refe{theo2}  $.

\medskip

Let us now say a few words about the simple idea on which Theorem
 $\refe{theo1} $ is based. The proof  uses Mourre's theory and
relies on two remarks. The first one is roughly the following: assume that,
 for $  \lambda \gg 1  $, we can find $ f_{\lambda} \in C_0^{\infty} (\Ra) $ and some
self-adjoint operator $ A $ such that the (formal) commutator $ i
[H,A] $ has a bounded closure $ i [H,A]^0 $ on $ D (H) $  and
\begin{eqnarray}
 f_{\lambda} (H) i [H,A]^0 f_{\lambda} (H) \geq \lambda
f^2_{\lambda} (H) \label{Mourreroughly}
\end{eqnarray}
 with $ f_{\lambda} = 1 $ on
$ ( \lambda -  \delta_{\lambda} , \lambda + \delta_{\lambda} ) $.
Then, one has
$$ || \scal{A}^{-s}(H-\lambda \pm i 0)^{-1}
\scal{A}^{-s} || = {\mathcal O}(\delta_{\lambda}^{-1}) .$$ This
essentially follows from the techniques of \cite{Mour1} (thought
our assumptions on $A $ and $H$ won't fit the framework of
\cite{Mour1}) 
and is the purpose of the next section. We emphasize that, instead of $
\reff{Mourreroughly} $, a Mourre estimate usually looks like
\begin{eqnarray}
 E_{I(\lambda)} (H) i [H,A]^0 E_{I(\lambda)} (H) \geq 2 \lambda E_{I(\lambda)} (H) +
E_{I(\lambda)} (H) K_{\lambda} E_{I(\lambda)} (H)  \label{Mourrereste}
\end{eqnarray}
 with $
E_{I(\lambda)} (H) $ the spectral projector of $ H $ on some
interval $ I(\lambda) \ni \lambda $, and $ K_{\lambda} $ a compact
operator. As explained in \cite{Mour1}, $ \reff{Mourrereste} $ implies $
\reff{Mourreroughly} $ provided $ f_{\lambda} $ is supported away
from the point spectrum of $ H $ and $ \delta_{\lambda} $ is small
enough, since $ f_{\lambda}(H)K_{\lambda} \rightarrow 0 $ as $
\delta_{\lambda} \rightarrow 0 $. But {\it we don't have any
control on} $ \delta_{\lambda} $ in general and here comes our
second remark. If one already knows some a priori estimates on $
(H - \lambda \pm i 0)^{-1} $, we can hope to control $
\delta_{\lambda} $ from below by mean of the following easy lemma
which links explicitly the size of the support of the function,
i.e. $ \delta_{\lambda} $, to estimates on the resolvent.

\begin{lemm} \label{easytrick} Let $ (L,D(L)) $ be a self-adjoint operator on
  a Hilbert space $ {\mathcal H} $ and $ J $ an interval. Assume that, for some
  bounded operator $ K $,
\begin{eqnarray}
 \sup_{\lambda \in J, \ 0 <\epsilon < 1} 
\left| \left| K^* ( L - \lambda \pm i \epsilon)^{-1} K \right|
\right| < \infty  . \label{Lebesgueuniforme} 
\end{eqnarray}
 Then, for all
$ f \in C_0^{\infty}(J) $, one has
$$ || f ( L ) K || \leq \pi^{-1/2} |J|^{1/2} || f ||_{\infty} \ \sup_{  \lambda \in J}
\left| \left| K^* ( L - \lambda \pm i 0)^{-1} K \right|
\right|^{1/2} ,
$$ with
$ | J | $ the Lebesgue measure of $ J $, provided the right hand side is
well defined.
\end{lemm}

\noindent {\it Proof.}  This is a direct consequence of the
Spectral Theorem which shows  that, for all $ \varphi  \in
{\mathcal H} $,
$$ || f (L) K \varphi ||^2 = (2 i \pi)^{-1}
\lim_{\epsilon \downarrow 0} \int_J  | f (E)|^2 \left( \big( (L -
E - i \epsilon)^{-1} - (L- E + i \epsilon)^{-1} \big) K \varphi  ,
K \varphi  \right) \ \mbox{d}E. \qquad \Box $$

\medskip

\noindent {\bf Remark.} If $ L = H $ and $ J \Subset ((n-1)^2/4,+\infty)  $,
 the condition $ \reff{Lebesgueuniforme} $ is
known to hold by \cite{CaVo1,FrHi1}, choosing
for instance $ K = \scal{r}^{-s} $ with $ s > 1/2  $.

\medskip

 We shall apply this strategy, i.e. deduce $ \reff{Mourreroughly} $ from an
 estimate of the type  $ \reff{Mourrereste} $ using the above trick with the a
priori estimates of Cardoso-Vodev  proved in \cite{CaVo1}. The
conjugate operator $ A $ (which will actually depend on $ \lambda
$) is essentially the one constructed by Froese-Hislop in \cite{FrHi1}.

We note in passing that we actually prove a stronger result than Theorem $
\refe{theo1} $, namely a Mourre estimate (see Theorem $ \refe{checkconditions2} $) which
implies Theorem $ \refe{theo1} $. Thus, using the techniques of \cite{Mour2}, we could 
also get other propagation estimates involving "incoming" or "outgoing"
spectral cutoffs.

This method is rather general and could certainly be adapted to
other settings
  than the asymptotically hyperbolic
one. For instance, we could consider manifolds with Euclidean ends
or both asymptotically hyperbolic and Euclidean ends, using the
standard generator of dilations $ r D_r + D_r r $ (cut off near
infinity) as a conjugate operator in Euclidean ends, as in
\cite{FrHi1}.


\medskip

The organization of the paper is the following. In Section $
\refe{Mourrerevise} $, we review Mourre's theory with a class of 
operators adapted to our purpose  and give a rather
explicit dependence of the estimates with respect to the different
parameters.  We point out that some of our technical assumptions
on $ A $ and $ H $ will not be the same as those of \cite{Mour1}.
For this reason and also to take the parameters into account, we
need to provide some details. In Section $ \refe{opconjugue} $, we
review the construction of the conjugate operator $A $ introduced
in \cite{FrHi1}. For the same reasons as for Section $
\refe{Mourrerevise} $, we cannot use directly the results of
\cite{FrHi1} and we need again to review some proofs. We also give
a pseudo-differential approximation for $ A $. In Section $
\refe{application} $, we prove Theorems $ \refe{theo1} $ and $
\refe{theo2} $.

\section{Mourre's theory } \label{Mourrerevise}
\setcounter{equation}{0}
\subsection{Algebraic results} \label{conditions}
In what follows, $ (H,D(H)) $ and $ (A,D(A)) $ are self-adjoint
operators on a Hilbert space $ {\mathcal H} $  that will
eventually satisfy the assumptions {\it (a)}, {\it (b)} and {\it
(c)} below. These assumptions are slightly different from the ones
used in \cite{Mour1} but, taking into account some minor
modifications, they allow to follow the original proof of Mourre
to get estimates on $ \scal{A}^{-s}(H-\lambda \pm i 0)^{-1}
\scal{A}^{-s} $.
 In this subsection, we record results allowing
 to justify the algebraic manipulations needed for that purpose.
Differential  inequalities and related estimates are given in Subsection
$ \refe{differentialinequalities} $.

\smallskip

 \noindent {\it (a) Assumptions on domains:} there exists a subspace $
{\mathcal D} \subset D (H) \cap D (A)  $ dense in  $ {\mathcal H}
$, such that
\begin{eqnarray}
{\mathcal D} \ \mbox{ is a core for } A, \label{coreA}
\end{eqnarray}
i.e. is dense in $ D(A) $ equipped with the graph norm. We also
assume  the existence of a sequence $ \zeta_n $ of bounded
operators satisfying, for all $n \in \Na$,
\begin{eqnarray}
\zeta_n D (H) \subset D (H) , \qquad \zeta_n D (A) \subset D (A) ,
\label{stabiliteAH} \\
 \zeta_n (H-z)^{-1} {\mathcal D} \subset {\mathcal D}, \qquad \forall \ z
\notin \mbox{spec} (H), \label{regularity} \\
 \zeta_n g (H) {\mathcal H} \subset {\mathcal D}, \qquad \forall \ g \in
C_0^{\infty}(\Ra), \label{regularityinf}
\end{eqnarray}
and furthermore, as $ n \rightarrow \infty $,
\begin{eqnarray}
 \zeta_n \varphi & \rightarrow & \varphi , \qquad \forall \ \varphi \in
{\mathcal H} ,
\label{strong} \\
 A   \zeta_n  \varphi & \rightarrow & A \varphi , \qquad \forall \ \varphi \in
 D (A),
\label{strongA} \\
 H  \zeta_n  \varphi & \rightarrow & H \varphi  , \qquad \forall \ \varphi \in
 D (H) .
\label{strongH}
\end{eqnarray}
The last condition regarding the domains is the following
important one
\begin{eqnarray}
(H-z)^{-1}D(A) \subset D (A), \qquad \forall \ z \notin
\mbox{spec} (H) . \label{stable}
\end{eqnarray}

\medskip

\noindent {\bf Remark.} When $A$ and $H $ are pseudo-differential
operators on manifolds, most of these conditions are easily
verified. The hardest is to check $ \reff{stable} $. We point out
that sufficient conditions ensuring $ \reff{stable} $ are given in
\cite{Mour1} (see also \cite{ABG1,GeGe1}), namely conditions on $ e^{itA} $, 
but they don't seem to be satisfied by the operators considered in
Section $ \refe{opconjugue} $. We thus rather set \reff{stable} as
an assumption in this part; in the next section, the explicit
forms of $ A $ and $ H $ will allow us to check it directly 
(see Proposition $ \refe{checkconditions} $).

\medskip

Note also the following easy result.
\begin{lemm} \label{commoncore} Conditions $ \reff{stabiliteAH} $, $
  \reff{strong} $, $ \reff{strongA} $ and $ \reff{strongH}
$  imply that $ A \zeta_n (A+i)^{-1} $, $H \zeta_n (H+i)^{-1} $
are bounded operators on $ {\mathcal H} $, uniformly with respect to $ n $. In
addition,  $ \reff{regularity} $ implies that  $ {\mathcal D} $ is a
core for $ H $.
\end{lemm}

\noindent {\it Proof.} We only consider   $ H $. For all $
\epsilon > 0 $, $ H (\epsilon H + i)^{-1} \zeta_n (H + i)^{-1} $
is bounded and converges strongly on $ {\mathcal H} $ as $
\epsilon \rightarrow 0 $, since $D (H) $ is stable by $ \zeta_n $.
This proves that  $ H \zeta_n (H +i)^{-1} $ is bounded, by uniform
boundedness principle. Then, by $ \reff{strongH} $, $ H \zeta_n
(H+i)^{-1} $ converges strongly on $ {\mathcal H} $ to $ H
(H+i)^{-1} $ and hence is uniformly bounded by the same principle.
Thus, if $
 \psi \in D (H) $ and $ {\mathcal D} \ni \varphi_n  \rightarrow
(H+i) \psi $ in $ {\mathcal H} $, then $ \psi_n := \zeta_n
(H+i)^{-1} \varphi_n $ is clearly  a sequence of $ {\mathcal D} $ such that
$ \psi_n \rightarrow \psi  $ and $ H \psi_n \rightarrow H \psi $
in $ {\mathcal H} $. \finpreuve

\medskip

\noindent {\it (b) Commutators assumptions.} There exists a
bounded operator $ [H,A]^0 $ from $ D (H) $ (equipped with the
graph norm) to $ {\mathcal H} $, and $ C_{H,A} > 0 $ such that,
for all $ \varphi,\psi \in {\mathcal D} $,
\begin{eqnarray}
(A \varphi , H \psi) - (H \varphi , A \psi) = \left( [H,A]^0
\varphi , \psi \right), \qquad \qquad \qquad  \label{comm1}  \\
\left| \left( A \varphi , i [H,A]^0 \psi \right) - \left( i
[H,A]^0 \varphi , A \psi \right) \right| \leq C_{H,A}|| \psi|| \
|| (H+i) \varphi || .  \label{comm2}
\end{eqnarray}

Note that we only require that $ \varphi , \psi \in {\mathcal D} $
in $ \reff{comm1}$ and $\reff{comm2} $ (instead of $ D (A) \cap D
(H) $ in the original paper \cite{Mour1}). Note also  that $ i
[H,A]^0 $ is automatically {\it symmetric} on $ {\mathcal D} $,
hence on $ D (H) $  by Lemma $ \refe{commoncore} $.

\medskip

We now state the main assumption.

\noindent {\it (c) Positive commutator estimate at $ \lambda \in
\Ra $.} There exists  $ \delta > 0 $ and $ f \in
C_0^{\infty}(\Ra,\Ra) $ with $ 0 \leq f \leq 1$,  such that,
\begin{eqnarray}
 f (E)  =
  \begin{cases}
    1 & \text{if} \ \ |E-\lambda| < 2 \delta , \\
    0 & \text{if} \ \  |E-\lambda| > 3 \delta ,
  \end{cases}  \nonumber 
\end{eqnarray}
and satisfying, for some $ \alpha > 0 $,
\begin{eqnarray}
f(H) i [H,A]^0 f (H) \geq \alpha f (H)^2 . \label{commpositif}
\end{eqnarray}
 Remark that $ \reff{commpositif} $ makes perfectly
sense, for $ f (H) i [H,A]^0 f (H) $ is bounded and self-adjoint
in view of the symmetry of $ i [H,A]^0 $ on $ D (H) $.

The main condition among {\it (a)}, {\it (b)} and {\it (c)} is the
{\it Mourre estimate} $ \reff{commpositif} $.  We include the
parameters $ \alpha $ and $ \delta $ to emphasize their important
roles in the estimates given in the next subsection.

\bigskip

We now record the main  algebraic tools needed to repeat Mourre's
strategy.

\begin{prop} \label{blambda} Assume that all the conditions $ \reff{coreA},
\cdots, \reff{comm1} $ but $ \reff{regularityinf}  $  hold. Then,
on $ D (A) $,
\begin{eqnarray}
 [(H-z)^{-1},A] = -
(H-z)^{-1}[H,A]^0 (H-z)^{-1}, \qquad z \notin \emph{spec} (H) .
\label{comm11}
\end{eqnarray}
Furthermore, $ (A \pm i \Lambda )^{-1} D (H) \subset D (H) $ for all  $
\Lambda \gg 1 $ and, by  setting $ A (\Lambda) = i \Lambda A (A + i
\Lambda)^{-1} $, we have
\begin{eqnarray}
 [H,A(\Lambda)] \varphi \rightarrow [H,A]^0  \varphi ,
 \qquad \Lambda \rightarrow \infty, \label{virial}
\end{eqnarray}
 in $ {\mathcal H} $, for all $ \varphi \in D (H) $.
\end{prop}

\noindent {\it Proof.} We apply $ \reff{comm1} $ to $ \varphi_n =
\zeta_n (H-z)^{-1} \tilde{\varphi} $ and $ \psi_n = \zeta_n
(H-\bar{z})^{-1} \tilde{\psi} $ with $ \tilde{\varphi},
\tilde{\psi} \in {\mathcal D} $. Since $ [H,A]^0 $ is bounded on $
D (H) $, $ \reff{strongH} $ implies that  $ [H,A]^0 \varphi_n \rightarrow
[H,A]^0 (H-z)^{-1} \tilde{\varphi} $. Furthermore, $ \zeta_n
(H-z)^{-1} \tilde{\varphi} \rightarrow (H-z)^{-1} \tilde{\varphi}
$ in $ D (A) $ by $ \reff{strongA} $ and $ \reff{stable} $ (the
same holds for $ \tilde{\psi} $) and hence
$$ \left((H-z)^{-1} \tilde{\varphi} , A \tilde{\psi}
\right) - \left(A \tilde{\varphi} , (H-\bar{z})^{-1} \tilde{\psi}
\right) = \left( [H,A]^0 (H-z)^{-1} \tilde{\varphi},
(H-\bar{z})^{-1} \tilde{\psi} \right) .
$$
Since $ {\mathcal D} $ is a core for $ A $, the above equality
actually holds for all $ \tilde{\varphi}, \tilde{\psi} \in D (A) $. This shows
 $ \reff{comm11} $. The proof of $ \reff{virial}
$ follows as in  \cite{Mour1}. Indeed  $ \reff{comm11} $ yields
\begin{eqnarray}
 [(H-z)^{-1},(A-Z)^{-1}] = - (A-Z)^{-1} (H-z)^{-1} [H,A]^0
(H-z)^{-1}(A-Z)^{-1} , \label{resol}
\end{eqnarray}
which implies that $ (H+i)^{-1}(A \pm i \Lambda)^{-1} = (A \pm i
\Lambda)^{-1} (H+i)^{-1}(1 + {\mathcal O}(\Lambda^{-1})) $, where
$ {\mathcal O}(\Lambda^{-1}) $ holds in the operator sense. This
clearly implies that $  (A \pm i \Lambda )^{-1} D (H) \subset D
(H) $ for $ \Lambda \gg 1  $ and that $ B (\Lambda):= (H+i) i
\Lambda (A + i \Lambda)^{-1}(H+i)^{-1} \rightarrow 1 $, in the
strong sense on $ {\mathcal H} $.  The latter leads to $
\reff{virial} $ since, on $ D (H) $,
\begin{eqnarray}
 [H,A(\Lambda)] = i \Lambda (A + i \Lambda)^{-1} [H,A]^0
(H+i)^{-1} B (\Lambda) (H+i) . \label{virialbis}
\end{eqnarray}
 The proof is complete. \finpreuve

\bigskip

The next proposition is important for several reasons. Firstly, it
will allow to justify the manipulation of some commutators and
secondly, it gives an explicit estimate for the norm of (the
closure of) $ [g(H),A](H+i)^{-1} $. It is also a key to the proof
of the useful Proposition $ \refe{adap} $ below. We include the
proof of Proposition $ \refe{derivecomm} $, essentially taken from
\cite{Mour1}, to convince the reader that our assumptions are
sufficient to get it.

\begin{prop} \label{derivecomm} Under the assumptions of
Proposition $ \refe{blambda} $, the following holds: for any
bounded Borel function $ g $ such that $ \int |t \hat{g}(t)|
\emph{d}t < \infty $, we have $ g (H) (  D (A) \cap D (H) ) \subset D (A)
$ and
$$ \left| \left| [g(H),A] \varphi \right| \right| \leq (2 \pi)^{-1} \int |t \hat{g}(t)|
\emph{d}t\ ||[H,A]^0 (H+i)^{-1}|| \ || (H+i) \varphi ||, \qquad
\forall \ \varphi \in D (A) \cap D (H) .
$$
\end{prop}
Before proving this proposition, we quote the following important
consequence.
\begin{prop} \label{adap} In addition to the assumptions of Proposition $ \refe{blambda} $, suppose that
 $ \reff{regularityinf} $ holds. Then, for any $ \varphi \in D (A) \cap D (H) $, there exists a
sequence $ \varphi_n  \in {\mathcal D} $ such that, as $ n
\rightarrow \infty $,
$$ \varphi_n \rightarrow \varphi, \qquad
 A \varphi_n \rightarrow A \varphi \qquad \mbox{and} \qquad H \varphi_n \rightarrow H \varphi . $$
In particular, $ \reff{comm1} $ and $ \reff{comm2} $ hold for all
$ \varphi, \psi \in D (A) \cap D (H) $.
\end{prop}

\noindent {\it Proof.}  We choose $ g \in C_0^{\infty}(\Ra)$,  $ g
= 1 $ near $ 0 $, and set $ \varphi_n = \zeta_n g_n (H) \varphi $,
with $ g_n (E) = g (E/n) $. It belongs to $ {\mathcal D} $ by $
\reff{regularityinf} $ and clearly converges to $ \varphi $ in $
{\mathcal H} $. Furthermore, $ (H+i) \zeta_n (H+i)^{-1} $
converges strongly on $ {\mathcal H} $ by $ \reff{strongH} $ and
this easily shows that $ H \varphi_n \rightarrow H \varphi $.
Regarding $ A \varphi_n $, we write
$$ A
\varphi_n = A \zeta_n (A + i)^{-1} g_n (H) (A+i) \varphi - A
\zeta_n (A+i)^{-1} [g_n(H),A] \varphi $$ where $ A \zeta_n
(A+i)^{-1} $ converges strongly on $ {\mathcal H} $ by $
\reff{strongA} $ and $ ||[A,g_n (H)] \varphi || \leq C n^{-1} ||
(H+i) \varphi || $ since $ \int |t \hat{g}_n(t)| \mbox{d}t =
{\mathcal O}(n^{-1}) $.  \finpreuve
\medskip

As a consequence of this proposition, we can define, for further
use, the form $ [[H,A]^0,A] $ by
\begin{eqnarray}
 \left( [[H,A]^0,A] \varphi, \psi \right) : = \left( A \varphi , i
[H,A]^0 \psi \right) - \left( i [H,A]^0 \varphi , A \psi \right) ,
\qquad \varphi , \psi \in D (A) \cap D (H) .
\label{commdoublefull}
\end{eqnarray}

\noindent {\it Proof of Proposition $ {\it \refe{derivecomm}}$.}
We first observe that, if $ \varphi \in D (A) \cap D (H)  $ and $
\Lambda \gg 1 $, then for all $ t $
$$ e^{itH} A (\Lambda) e^{-itH} \varphi = A (\Lambda) \varphi + i \int_0^t e^{isH}
[H,A(\Lambda)] e^{-isH} \varphi \ \mbox{d}s . $$ This can be
easily seen  by weakly differentiating both sides with respect to
$ t $, testing them against an arbitrary element of $  D (H) $. This equality
shows that, for any $ \psi \in {\mathcal H} $,
\begin{eqnarray}
 \left( [A(\Lambda),g(H)] \varphi , \psi \right) = \frac{i}{2
\pi} \int \hat{g}(t) \int_0^t \left( e^{-i(t-s)H} [H,A(\Lambda)]
(H+i)^{-1} e^{-isH} (H+i) \varphi , \psi \right)  \mbox{d}s
\mbox{d}t . \label{clef2}
\end{eqnarray}
By $ \reff{virialbis} $, $ [H,A(\Lambda)](H+i)^{-1} $ is uniformly
bounded, so  the modulus of right hand side is dominated by $ C ||
\psi || $, for some $ C $ independent of $ \Lambda $. In
particular, if $ \psi \in D (A) $,
$$ (g(H) \varphi , A \psi) = \lim_{\Lambda \rightarrow \infty} (g(H),A(-\Lambda) \psi) =
\lim_{\Lambda \rightarrow \infty} (g(H)A(\Lambda)\varphi, \psi) -
([g(H),A(\Lambda)]\varphi,\psi) $$ proves that $ |(g(H) \varphi ,
A \psi )| \leq C || \psi || $, with $ C $ independent of $ \psi
\in D (A) $. This implies that $ g (H) \varphi \in D (A^*) = D
(A) $. Then, letting $ \Lambda \rightarrow \infty $ in $
\reff{clef2} $  clearly leads to the estimate on $|| [g(H),A]
\varphi ||$ . \finpreuve

\medskip

We now quote a crucial result which is directly taken from
\cite{Mour1}.

\begin{prop} \label{algebre} Assume that $ B $ is a bounded operator on $ {\mathcal H} $.
 Then for any $ z \notin \Ra $ and any $
\varepsilon \in \Ra $ such that $ \emph{Im}(z) \varepsilon \geq 0
$, the operator $ H - z - i \varepsilon B^* B $ is a bounded
isomorphism from $ D (H) $ (with the graph norm) onto $ {\mathcal
H} $. If we set
$$ G_z (\varepsilon) = (H-z-i \varepsilon B^* B)^{-1} $$
we have, provided $  \emph{Im}(z) \varepsilon \geq 0 $ and $
\emph{Im}(z) \varepsilon_0 \geq 0 $,
\begin{eqnarray*}
 G_z (\varepsilon) -  G_z (\varepsilon_0) = G_z (\varepsilon)i
(\varepsilon - \varepsilon_0)B^*B G_z (\varepsilon_0) ,  \\  G_z
(\varepsilon)^* = G_{\bar{z}}(-\varepsilon), \qquad ||G_z
(\varepsilon) ||
 \leq |\emph{Im}(z)|^{-1},
 \end{eqnarray*}
 in the sense of bounded operators on $ {\mathcal H} $.
 Furthermore, if $ B^{\prime} $ and $ C $ are bounded operators,
 with $ C $ self-adjoint, and if $ \emph{Im}(z)\varepsilon > 0 $,
 then
\begin{eqnarray}
  B^{\prime *} B^{\prime} \leq B^* B \qquad \Rightarrow \qquad \left| \left|
B^{\prime} G_z (\varepsilon) C
  \right| \right|
  \leq |\varepsilon|^{-1/2} \left| \left| C G_z (\varepsilon) C  \right|
\right|^{1/2} . \nonumber
  \end{eqnarray}
\end{prop}

This result, which is one of the keys of the differential
inequality technique of Mourre, will of course be used with $ B^*
B = f (H) i [H,A]^0 f (H) $, but it doesn't depend on any of the
assumptions quoted in the beginning of this section. We refer to
\cite{Mour1} for the proof and rather put emphasize on the
following result.

\begin{prop} Assume that all the conditions from $ \reff{coreA} $
to $ \reff{commpositif} $ hold and define $ G_z (\varepsilon) $ as
above with $ B^* B = f (H) i [H,A]^0 f (H) $. Then $ G_z
(\varepsilon) D (A) \subset D (A) \cap D (H) $.
\end{prop}

\noindent {\it Proof.} It suffices to show that $ G_z
(\varepsilon) \varphi  $ belongs to $ D (A) $ for any $ \varphi
\in D (A) $. As in the proof of Proposition $ \refe{derivecomm} $,
this is implied by the fact that $ \sup_{\Lambda \geq \Lambda_0}||
[ G_z (\varepsilon), A (\Lambda) ] || < \infty $, for $ \Lambda_0
$ large enough. To prove this, we remark that
$$ [A(\Lambda),G_z (\varepsilon)] = G_z (\varepsilon) [H,A(\Lambda)] G_z (\varepsilon) - i \varepsilon
G_z (\varepsilon) [B^*B,A(\Lambda)] G_z (\varepsilon) $$ where the
first term of the right hand side is uniformly bounded by $
\reff{virial} $ and the uniform boundedness principle. We are thus
left with the study of the second term for which we observe that
\begin{eqnarray}
 \left( [(A+i\Lambda)^{-1},B^*B] \psi_1 , \psi_2 \right)  & = &
\left( i [H,A]^0 f (H) \psi_1 , [f(H),(A-i\Lambda)^{-1}] \psi_2
\right) + \nonumber \qquad \\
& & \left( [(A+i\Lambda)^{-1},f(H)] \psi_1 , i [H,A]^0 f (H)
\psi_2 \right) + \nonumber \\ & &  \left( [(A+i\Lambda)^{-1}, i
[H,A]^0] f (H) \psi_1 , f (H) \psi_2 \right) , \nonumber
 \label{lastline}
\end{eqnarray}
for all $ \psi_1 , \psi_2 \in D (H) $. Since $ A (\Lambda) = i
\Lambda + \Lambda^2 (A + i \Lambda)^{-1} $, multiplying this
equality by $ \Lambda^2 $ allows to replace $ (A \pm i
\Lambda)^{-1} $ by $ A(\pm \Lambda) $. By  $ \reff{virialbis} $
and $
 \reff{clef2} $ , $ [f(H),A(\pm \Lambda)](H+i)^{-1} $ is
uniformly bounded which reduces the proof of the proposition to
the study of $ \left( [A(\Lambda), i [H,A]^0] f (H) \psi_1 , f (H)
\psi_2 \right) $. To that end, we note that, if $ \tilde{\psi}_1 ,
\tilde{\psi}_2 $ belong to $  D (H) $, then $ (
[A(\Lambda),i[H,A]^0] \tilde{\psi}_1 , \tilde{\psi}_2 ) $
 can be written
\begin{eqnarray}
\Lambda^2 \left( A (A + i \Lambda)^{-1} \tilde{\psi}_1 , i [H,A]^0
(A-i\Lambda)^{-1} \tilde{\psi}_2 \right) - \Lambda^2 \left( i
[H,A]^0 (A + i \Lambda)^{-1} \tilde{\psi}_1 , A (A - i
\Lambda)^{-1} \tilde{\psi}_2 \right). \nonumber
\end{eqnarray}
Using  $ \reff{comm2} $ and Proposition $ \refe{adap} $, combined
with the fact that $ \Lambda (H + i) (A \pm i \Lambda)^{-1}
(H+i)^{-1} $ is uniformly bounded (see the proof of Proposition $
\refe{blambda}$), we obtain the existence of $ C > 0 $ such that
$$ \left| \left( [A(\Lambda),i[H,A]^0] \tilde{\psi}_1 , \tilde{\psi}_2 \right) \right| \leq C
|| (H+i)\tilde{\psi}_1 || \ || \tilde{\psi}_2 ||, \qquad
\tilde{\psi}_1 , \tilde{\psi}_2 \in  D (H)
$$
for $ \Lambda \gg 1 $. The conclusion follows.
\finpreuve

Note that we have chosen to include this proof, thought it is essentially the one of
 \cite{Mour1}, since our assumptions on $ A $ are not the same as those of \cite{Mour1}.

\subsection{The limiting absorption principle} \label{differentialinequalities}
In this part, we repeat the method of differential inequalities of
Mourre \cite{Mour1} to get estimates on the boundary values of $
(H-z)^{-1} $. Our main goal is an explicit control of the
different estimates in terms of the parameters, namely $ A,H, f,
\lambda , \alpha , \delta $ and $ C_{H,A} $ (see $ \reff{comm2}
$). As we shall see, the following quantities will play a great
role
\begin{eqnarray}
N_{[H,A]} & : = & \left| \left| [H,A]^0 (H+i)^{-1} \right| \right| \\
S_{H,A}^{f,\alpha} & : = & \left( 1 + \alpha^{-1} \left| \left|
[H,A]^0 f (H) \right|
 \right|\right)^2 \\
\Delta_f & : = & (2 \pi)^{-1} \int_{\Ra} |t \hat{f}(t)| \
\mbox{d}t .
\end{eqnarray}

We assume that all the conditions from $ \reff{coreA} $ to $
\reff{commpositif} $ hold and that $ G_z (\varepsilon) $ is
defined by Proposition $ \refe{algebre} $ with $ B^* B = f (H) i
[H,A]^0 f (H) $.

\medskip

As a direct consequence of Proposition $ \refe{algebre} $, we
first get the  estimate
\begin{eqnarray}
\left| \left| f(H) (H+i)^k  G_z (\varepsilon) w(A) \right| \right|
\leq (1 + |\lambda| + 3 \delta )^k
\alpha^{-1/2}|\varepsilon|^{-1/2} \left| \left| w(A) G_z
(\varepsilon) w (A) \right| \right|^{1/2} \nonumber
\end{eqnarray}
which holds for any bounded and real valued Borel function $ w $.
We also obtain immediately
\begin{eqnarray}
 \left| \left| f (H) G_z (\varepsilon) f (H)  \right| \right|  \leq
\alpha^{-1}|\varepsilon|^{-1} . \label{unsurepsilon}
\end{eqnarray}
On the other hand, by the resolvent identity given in Proposition
$ \refe{algebre} $, we see that
$$  G_z (\varepsilon) f (H) =  G_z (0) \left( f (H) - \varepsilon f (H)  [H,A]^0 f (H)
G_z (\varepsilon) f (H) \right) $$ where the bracket is uniformly
bounded with respect to $ \varepsilon $ by $ \reff{unsurepsilon} $
and we obtain
\begin{eqnarray}
 \left| \left| (H+i)^k (1-f)(H) G_z (\varepsilon) f (H) \right| \right| \leq \sup_{|E-\lambda| \geq 2 \delta}
 \frac{|E+i|^k}{|E-z| }
 \left( 1 + \alpha^{-1} \left| \left| [H,A]^0 f (H) \right|
 \right|
 \right) \label{zeroepsilon}  ,
\end{eqnarray}
for $ k = 0,1 $. Here we used the fact that $ f(H) [H,A]^0 $ has a
bounded closure whose norm equals $ || [H,A]^0 f(H) || $. Another
application of the resolvent identity also gives
\begin{eqnarray}
 G_z (\varepsilon) (1-f)(H) = G_z (0) \left( (1 -f) (H) -
\varepsilon f (H)  [H,A]^0 f (H) G_z (\varepsilon) (1-f) (H)
\right) \label{zeroepsiloncarre}
\end{eqnarray}
 in which $ f (H) G_z (\varepsilon) (1-f) (H)  $ can be
 estimated (independently of $ \varepsilon $) using  $ \reff{zeroepsilon}
 $.

 Summing up, all this leads to
\begin{prop} \label{apriori4} Assume that  $\lambda , \delta , \alpha $ satisfy condition {\it (c)}
of Subsection $ \refe{conditions} $ and that
\begin{eqnarray}
\varepsilon \emph{Im} \ z > 0, \qquad | \emph{Re} \ z  - \lambda|
\leq \delta, \qquad \delta \leq \alpha \qquad \mbox{and} \qquad
|\varepsilon| \leq \delta \alpha^{-1}. \label{premierchoix}
\end{eqnarray}
Then, for $ k = 0,1 $ and all bounded Borel function $w$ such that
$ ||w||_{\infty} \leq 1 $, we have
\begin{eqnarray}
\left| \left| (H+i)^k (1-f)(H) G_z (\varepsilon) \right| \right| &
\leq &
(1+|\lambda|+2 \delta)^k \delta^{-1} \left( 1 + S_{H,A}^{f,\alpha} \right) , \label{Czero} \\
\left| \left| (H+i)^k f(H)   G_z (\varepsilon) w(A) \right|
\right| & \leq & (1 + |\lambda| + 3 \delta )^k
\alpha^{-1/2}|\varepsilon|^{-1/2} \left| \left| w(A) G_z
(\varepsilon) w (A) \right| \right|^{1/2} , \label{Chalf}
\\
 \left| \left| w(A) G_z
(\varepsilon) w(A)  \right| \right| & \leq & \alpha^{-1}
|\varepsilon|^{-1} \left(2 +  S_{H,A}^{f,\alpha} \right) .
\label{departdiff}
\end{eqnarray}
\end{prop}

Note that the right hand side of $ \reff{Czero} $ is independent
of $ \varepsilon $. Note that we also get estimates  on $ G_z (
\varepsilon ) (1-f)(H) $ and $ w (A) G_z (\varepsilon) f (H) $ for
free, by taking the adjoints, since $ G_z ( \varepsilon )^* =
G_{\bar{z}}(- \varepsilon) $.

\bigskip

We then need to get an estimate on $ d G_z (\varepsilon) / d
\varepsilon $. To that end, we simply  repeat the proof of Mourre
\cite{Mour1}, observing that the algebraic manipulations are valid
in our context thanks to the results of  Subsection $
\refe{algebre} $. In the sense of quadratic forms on $ D (A) $,
using in particular  $ [[H,A]^0,A] $ defined by $
\reff{commdoublefull} $, we thus obtain
\begin{eqnarray}
 \frac{d G_z (\varepsilon)}{d \varepsilon} & = & G_z
(\varepsilon) (1-f)(H)[H,A]^0 f (H) G_z (\varepsilon) + G_z
(\varepsilon) [H,A]^0 (1 - f)(H) G_z (\varepsilon) -  \nonumber
\\ & &
\varepsilon \left\{  G_z (\varepsilon) f (H)  [H,A]^0 [f(H),A] G_z
(\varepsilon) + G_z (\varepsilon) [f (H),A]  [H,A]^0 f(H) G_z
(\varepsilon) \right.  \nonumber \\ &  &  \left. + \  G_z
(\varepsilon) f (H) [ [H,A]^0,A] f(H) G_z (\varepsilon) \right\} +
  G_z
(\varepsilon) A  - A G_z (\varepsilon)  . \label{deriveeforme}
\end{eqnarray}
Let us set $ F_z (\varepsilon) : = w (A) G_z (\varepsilon) w (A)
$. By Proposition $ \refe{apriori4} $, $ \reff{deriveeforme} $
leads to the differential inequality
\begin{eqnarray}
\left| \left| w(A) \frac{ d G_z (\varepsilon)}{ d \varepsilon } w
(A) \right| \right| &  \leq  &  C_1 || F_z (\varepsilon) || +
C_{1/2} |\varepsilon|^{-1/2} || F_z (\varepsilon) ||^{1/2}  + C_0
\nonumber \\
& & \ \ + \ 2 || A w (A) || \left( \alpha^{-1/2}
|\varepsilon|^{-1/2} || F_z (\varepsilon) ||^{1/2} + \delta^{-1}
\left( 1 + S_{H,A}^{f,\alpha} \right) \right) \label{diff1}
\end{eqnarray}
where, by
Proposition $ \refe{derivecomm} $, the constants  $ C_0 , C_{1/2}$
and $ C_1 $ can be chosen as follows \label{constantespages}
\begin{eqnarray*}
C_0 & = &  \delta^{-2}  (1 + |\lambda| + 2 \delta )
 \left( 1 + S_{H,A}^{f,\alpha} \right)^2 N_{[H,A]} , \\
C_{1/2} & = & 2 \alpha^{-1/2} \delta^{-1} (1+|\lambda|+3 \delta)
S_{H,A}^{f,\alpha} N_{[H,A]}  \left( 1 + \delta \alpha^{-1}
\Delta_f N_{[H,A]} (1 + |\lambda| + 3 \delta)  \right),
\\
 C_1 & = & \alpha^{-1} (1+ | \lambda | + 3
\delta) \left(
 C_{H,A}  + 2 \Delta_f N_{[H,A]}^2 (1 + |\lambda | + 3 \delta )
 \right)   .
\end{eqnarray*}

 The second line of $ \reff{diff1} $
suggests that $ A w (A) $ must be bounded. Of course,  this holds
if $ w (E) = \scal{E}^{-1} $ (which was the original choice of
weight in \cite{Mour1}) however a trick of Mourre, which is
reproduced in \cite{PSS1},  allows to consider
$$ w (E) = \scal{E}_{\varepsilon}^{-s} := \scal{E}^{-s}
\scal{\varepsilon E}^{s-1}, \qquad 1/2 < s \leq 1 . $$ It is
indeed not hard to check that the following inequality holds for
all $ \varepsilon \ne 0 $ and $ E \in \Ra $
$$ \left| \frac{\partial}{\partial \varepsilon}
\scal{E}^{-s}_{\varepsilon}\right| = (1-s)
\scal{E}^{-s}_{\varepsilon}  \frac{ | \varepsilon | E^2}{1 +
\varepsilon^2 E^2}  \leq (1-s)|\varepsilon|^{s-1},
$$
and this implies that
\begin{eqnarray}
\left| \left| d \scal{A}^{-s}_{\varepsilon} / d \varepsilon  G_z
(\varepsilon) \scal{A}^{-s}_{\varepsilon} \right| \right|  \leq
 (1-s) |\varepsilon|^{s-1} \left( \alpha^{-1/2} |\varepsilon|^{-1/2} || F_z (\varepsilon) ||^{1/2}
  + \delta^{-1} \left( 1 + S_{H,A}^{f,\alpha} \right) \right) . \label{Leibnitz}
\end{eqnarray}
Using $ \reff{diff1} $,  $ \reff{Leibnitz} $ and the fact that $
\scal{E}^{-s}_{\varepsilon} \scal{E} \leq |\varepsilon|^{s-1} $
for $ 0 < |\varepsilon| \leq 1 $, we get the final differential
inequality
\begin{eqnarray}
\left| \left|  d F_z (\varepsilon) /d \varepsilon \right| \right|
& \leq &  C_1 || F_z (\varepsilon) || + C_{1/2}
|\varepsilon|^{-1/2} || F_z (\varepsilon) ||^{1/2}  + C_0
\nonumber \\
& & \ \  + \ 2 (2 - s) |\varepsilon|^{s-1} \left( 2 \alpha^{-1/2}
|\varepsilon|^{-1/2} || F_z (\varepsilon) ||^{1/2} + \delta^{-1}
\left( 1 + S_{H,A}^{f,\alpha} \right) \right)
 \label{diffs}
\end{eqnarray}
which is valid if  $ 0 < |\varepsilon| \leq 1 $ and if $
\reff{premierchoix} $ holds.

 Starting from  $
\reff{departdiff} $ and using $ \reff{diffs} $, a finite number of
integrations leads to a uniform bound on $ || F_z ( \varepsilon)
|| $ for $ 0 < |\varepsilon| \leq \min (1 , \delta \alpha^{-1} ) $
and thus on $ || F_z (0) || $. Such estimates depend of course on
$ A , H , f , \alpha , \lambda , \delta , C_0, C_{1/2} $ and $ C_1
$, but there is no reasonable way to express this dependence in
general. We thus rather consider a particular case in the  following
theorem, which lightens the role of $ \alpha , \lambda , \delta $.

\begin{theo} \label{theoclef} Consider families of operators $ H_{\nu} , A_{\nu} $, of
numbers $ \lambda_{\nu} , \alpha_{\nu} , \delta_{\nu} $ and of
functions $ f_{\nu} $ satisfying conditions {\it (a),(b),(c)} for
all $ \nu $ describing some set $ \Sigma $. Denote by $ C_{0,\nu},
C_{1/2,\nu} $ and $ C_{1,\nu} $ the corresponding constants
defined on page $ \pageref{constantespages} $.
 Assume that $ \varepsilon_{\nu}:= \delta_{\nu} \alpha_{\nu}^{-1} \leq 1 $ and
that there exists  $ C > 0 $ such that, for all $ \nu \in \Sigma
$,
\begin{eqnarray}
 C_{0,\nu} \leq C \varepsilon_{\nu}^{-1} \delta_{\nu}^{-1}, \! \! \qquad
C_{1/2,\nu} \leq C \varepsilon_{\nu}^{-1/2} \delta_{\nu}^{-1/2},
\! \! \qquad C_{1,\nu} \leq C \varepsilon_{\nu}^{-1}, \! \! \!
\qquad  || [H_{\nu},A_{\nu}]^0 f_{\nu}(H_{\nu})|| \leq C
\alpha_{\nu}
\label{constantesutiles}
\end{eqnarray}
with $ f_{\nu} $ of the form $ f_{\nu}(E) = f
((E-\lambda_{\nu})/\delta_{\nu}) $, for some fixed $ f \in
C_0^{\infty}(\Ra) $.
 Then, for all $ 1/2 < s \leq 1 $, there exists $ C_s
> 0 $ such that, for all $ \nu \in \Sigma $,
\begin{eqnarray}
 || \scal{A_{\nu}}^{-s}(H_{\nu}-z)^{-1} \scal{A_{\nu}}^{-s} ||
\leq C_s \delta_{\nu}^{-1} , \label{bornefinale}
\end{eqnarray}
 provided $ | \emph{Re} \ z -
\lambda_{\nu} | \leq \delta_{\nu} $. Furthermore, for any $ \mu
\in (\lambda_{\nu} - \delta_{\nu}, \lambda_{\nu} + \delta_{\nu})
$, the limits $$ \scal{A_{\nu}}^{-s} (H_{\nu}-\mu \pm i 0)^{-1}
\scal{A_{\nu}}^{-s} := \lim_{\varepsilon \rightarrow 0^+}
\scal{A_{\nu}}^{-s}(H_{\nu}-\mu \pm i \varepsilon)^{-1}
\scal{A_{\nu}}^{-s} $$ exist and are continuous, with respect to $
\mu $, in the operator topology.
\end{theo}

In practice, the conditions  $ \reff{constantesutiles} $ can be
checked using the explicit forms of $ C_0 , C_{1/2} $ and $ C_1 $
given on page $ \pageref{constantespages} $. We shall use this
extensively in the next section.

\medskip

\noindent {\it Proof.} We only consider the case where $
\varepsilon \in (0,\varepsilon_{\nu} ] $, i.e. the situation where
$ \mbox{Im } z$ is positive, since the one of $  \varepsilon \in
[-\varepsilon_{\nu} , 0 ) $ is similar. By the assumption on $ ||
[H_{\nu},A_{\nu}]^0 f_{\nu}(H_{\nu})|| $, the estimate $
\reff{departdiff} $ takes the form  $ || F_z (\varepsilon) || \leq
C \alpha_{\nu}^{-1} \varepsilon^{-1} $, thus $ \reff{diffs} $
implies that
$$ || F_{z}(\varepsilon) - F_z (\varepsilon_{\nu}) || \leq C_s \left( \delta^{-1}_{\nu}
+ \delta_{\nu}^{-1} \log (\varepsilon_{\nu}/\varepsilon) +
\alpha_{\nu}^{-1} \varepsilon^{s-1} \right), \qquad \forall \ \nu
\in  \Sigma  ,
$$
if $ 1/2 < s < 1 $. If $ s = 1 $, the term $ \varepsilon^{s-1} $
must be replaced by $ \log (\varepsilon_{\nu}/\varepsilon) $ which
can be absorbed by the second term of the bracket, for we assume
that $ \alpha_{\nu}^{-1} \leq \delta_{\nu}^{-1} $. Since $ || F_z
(\varepsilon_{\nu})|| \leq C \delta_{\nu}^{-1} $, a finite number
of iterations of Lemma $ \refe{recurrence} $ below completes the
proof of $ \reff{bornefinale} $. For the existence of the boundary
values of the resolvent, which are purely local, we refer to $
\cite{PSS1} $ (Theorem $8.1$). \finpreuve
\begin{lemm} \label{recurrence}  Let $ 0 \leq \sigma < 1 $ and
assume the existence of $ C $ such that, for all $ \nu \in \Sigma
$ and all $ \varepsilon\in (0 , \varepsilon_{\nu} ] $,
$$ || F_z (\varepsilon) || \leq C \left( \delta_{\nu}^{-1} + \delta_{\nu}^{-1}
 \log ( \varepsilon_{\nu} /
 \varepsilon ) +
\alpha_{\nu}^{-1}  \varepsilon^{-\sigma}   \right) . $$
 Then,
there exists $ C_{s,\sigma} $ such that, for all $ \nu \in \Sigma
$ and all $ \varepsilon \in ( 0 , \varepsilon_{\nu} ] $
$$ || F_z (\varepsilon) || \leq C_{s,\sigma} \left\{ \begin{array}{cc}
\delta_{\nu}^{-1} + \alpha_{\nu}^{-1} \varepsilon^{ s - 1/2 -
\sigma/2 }, & \mbox{if } \ s -1/2 < \sigma / 2 , \\
\delta_{\nu}^{-1} + \delta_{\nu}^{-1} \log ( \varepsilon_{\nu} /
\varepsilon ), & \mbox{if } \ s -1/2 = \sigma / 2, \\
\delta_{\nu}^{-1} , & \mbox{if } \ s-1/2 > \sigma/2 .
\end{array} \right.
$$
\end{lemm}

\noindent {\it Proof.} It simply follows from $ \reff{diffs} $ and
the fact that $ ||F_z (\varepsilon_{\nu})|| \leq C
\delta_{\nu}^{-1} $, by studying separately the three cases  and
using the trivial inequality
$$ \left( \delta_{\nu}^{-1} + \delta_{\nu}^{-1}
 \log ( \varepsilon_{\nu} /
 \varepsilon ) +
\alpha_{\nu}^{-1}  \varepsilon^{-\sigma} \right)^{1/2} \leq \
\delta_{\nu}^{-1/2} + \delta_{\nu}^{-1/2}
 \log^{1/2} ( \varepsilon_{\nu} /
 \varepsilon ) +
\alpha_{\nu}^{-1/2}  \varepsilon^{-\sigma/2} $$ to control the
terms involving $ || F_z (\varepsilon) ||^{1/2} $. \finpreuve

\section{Applications to asymptotically hyperbolic manifolds} \label{opconjugue}
\setcounter{equation}{0}
\subsection{The conjugate operator}
In this part, we recall the construction of the conjugate operator
defined by Froese-Hislop in \cite{FrHi1}. We emphasize that the
main ideas, namely the form of the conjugate operator and the
existence of a positive commutator estimate, are taken from
\cite{FrHi1}. However, since some of our assumptions (especially
{\it (a), (b)}  in subsection $ \refe{conditions} $) differ from
those of \cite{FrHi1} and since we need to control estimates with
respect to the spectral parameter, we will give a rather detailed
construction.

 Let $ \chi , \xi \in C^{\infty} (\Ra) $ be non
negative and non decreasing functions such that
$$  \chi (r) =  \begin{cases}
    0, & r \leq 1, \\
    1, & r \geq 2
  \end{cases} , \qquad  \xi (r) =  \begin{cases}
    0, & r \leq -1, \\
    1, & r \geq -\frac{1}{2}
  \end{cases} . $$
 By possibly replacing $ \chi $ and $ \xi $ by $ \chi^2 $ and $
\xi^2 $, we may assume that $ \chi^{1/2} $ and $ \xi^{1/2} $ are
smooth. For $ R
> r_0  $ and $ S > R $, we set $ \chi_R (r) = \chi (r/R) $ and $ \xi_S (r) = \xi (r/S)
$. Then, recalling that $ (\mu_k)_{k \geq 0} =
\mbox{spec}(\Delta_h) $ and setting $ \nu_k = (1 + \mu_k)^{1/2} $,
we define the sequence of smooth functions
$$ a_k (r) = (r + 2 S - \log \nu_k)\chi_R (r) \xi_S (r - \log \nu_k) .  $$
They are real valued and it is easy to check that their derivatives satisfy, for all $ j
\geq 1 $ and $ \ k \in \Na $,
\begin{eqnarray}
 | | a_k^{(j)} ||_{\infty} \leq C_j S^{1-j} , \qquad | |
a_k a_k^{(j+1)} | |_{\infty} \leq C_j  , \label{uniff}
\end{eqnarray}
uniformly with respect to  $ R > S >  r_0  $. Further on,  $ R $
and $ S $ will depend on the large spectral parameter $ \lambda $
but till then we won't mention the dependence of $ a_k $ (nor of
the related objects) on $ R,S $.

According to the results recalled in Appendix $ \refe{dimension1}
$,  there exists, for each $ k $,  a strongly continuous unitary
group $ e^{itA_k} $ on $ L^2 (\Ra) $ whose self-adjoint generator
$ A_k $ is
\begin{eqnarray}
 A_k =   a_k D_r  - i a^{\prime}_k  / 2  , \label{expressionak}
\end{eqnarray}
  i.e. a
self-adjoint realization of the r.h.s. Furthermore, we can
consider $ e^{itA_k} $ as a group on $ L^2 (I) $, since $
e^{itA_k} $ acts as the identity on functions supported in $
(-\infty , R) $ hence maps functions supported in $ I $ into
functions supported in $ I $ (see Appendix $ \refe{dimension1} $).
Therefore, using the notation $ \reff{coefFourier} $ for $
\varphi_k $, the linear map
\begin{eqnarray}
 \varphi \mapsto \sum_{k \geq 0} e^{itA_k} \varphi_k \otimes
\psi_k \label{sousunitaire}
\end{eqnarray}
 clearly defines a strongly continuous unitary group on $
L^2 (I) \otimes L^2 (Y, d \mbox{Vol}_h) $. The pull back on $ L^2 ({\mathcal M} \setminus
{\mathcal K}) $ of the operator $
\reff{sousunitaire} $, extended as the identity on $ L^2 ({\mathcal
K}) $, is also a strongly continuous unitary group on $ L^2
({\mathcal M}) $ which we denote by $ U (t) $. (Here again we omit
the $R,S$ dependence in the notation). Using Stone's Theorem \cite{RS1}, we
can state the
\begin{defi} \label{definition} We call $ A $ the self-adjoint generator of $ U (t) $.
In particular, its domain is  $$ D (A) = \{ \varphi \in L^2
({\mathcal M}) \ | \ U (t) \varphi \ \mbox{is strongly
differentiable at }  t =0 \},
$$
and $ A \varphi = i^{-1} d U (t) \varphi / dt_{ | t=0} $ for all $
\varphi \in D (A) $.
\end{defi}

\noindent {\bf Remark.} Note that this definition clearly implies
that $ L^2 ({\mathcal K}) \subset D (A) $ and that $ A_{| L^2
({\mathcal K})} \equiv 0 $.

\medskip

Now we choose  a sequence of functions  $   \zeta_n  \in
C_c^{\infty}({\mathcal M}) $ such that $ \zeta_n \rightarrow 1 $
strongly on $ L^2 ({\mathcal M}) $. More precisely, we choose $
\zeta_n $ of the form $ \zeta_n = \zeta (2^{-n} r ) $ for some $
\zeta \in C_0^{\infty}(\Ra) $ such that $ \zeta = 1 $ on a large
enough compact set (containing $0$) to ensure that $ \zeta_n = 1 $
near $ \overline{\mathcal K} $.

\begin{prop} \label{cbinf1} i) For all $ n $, $ \zeta_n D (A) \subset D (A)
$. \newline ii) For all $ \varphi \in D (A) $, $ A \zeta_n \varphi
\rightarrow A \varphi $ as $ n \rightarrow \infty $. \newline iii)
$  C_B^{\infty}( {\mathcal M}) $ is a core for $ A $ and $ A
C_B^{\infty}({\mathcal M}) \subset C_0^{\infty}({\mathcal M}) $.
\end{prop}

\noindent {\it Proof.} In view of the remark above, we only have
to consider $ \varphi \in L^2 ({\mathcal M} \setminus {\mathcal
K}) $ (i.e. supported in $ {\mathcal M} \setminus {\mathcal K} $).
Furthermore, to simplify the notations, we shall denote
indifferently by $ \varphi $ an element of $ L^2 ({\mathcal M}
\setminus {\mathcal K}) $ and the corresponding element in $
L^2(I) \otimes L^2 (Y,d \mbox{Vol}_h) $ via $ \reff{unitary} $.

Let us first observe that, for all such $ \varphi ,
\tilde{\varphi} $, Parseval's identity yields
$$ \left| \left| \left( U (t) \varphi  - \varphi \right) / it - \tilde{ \varphi} \right| \right|^2
= \sum_k \left| \left| \left( e^{itA_k} \varphi_k - \varphi_k
\right) / it - \tilde{\varphi}_k  \right| \right|^2 .
$$
Thus, by dominated convergence, this easily implies that $ \varphi
\in D (A) $ if and only if $ \varphi_k \in D (A_k) $ for all $ k $
and $ \sum_k ||A_k \varphi_k ||^2 < \infty  $, in which case $ (A
\varphi)_k = A_k \varphi_k $ for all $k$. Combining this
characterization with $ \reff{cutoffloin} $, and using the fact
that $ (a_k \zeta_n^{\prime})(r) = 2^{-n} a_k (r)
\zeta^{\prime}(2^{-n}r ) $ is uniformly bounded with respect to
$k,n \in \Na $ on $ I $, which is due to the fact that $ a_k (r)/r
$ is bounded with respect $r$ and $k$, we get {\it i)}. This also
shows that
$$  || A \zeta_n \varphi - \zeta_n A \varphi ||^2 = \sum_k || a_k \zeta_n^{\prime} \varphi_k ||^{2} $$
where the right hand side goes to $ 0 $ as  $ n \rightarrow \infty
$ by dominated convergence, and hence implies {\it ii)}. We now
prove {\it iii)}. Since $ A \varphi \equiv 0 $ for any function
supported outside $ \iota^{-1} ([R ,\infty)\times Y) $, and since
any element of $ D (A) $ can be approached by compactly supported
ones by {\it ii)}, it is clearly enough to show that for any $
\varphi \in D (A) $, compactly supported in $ \iota^{-1}
([R^{\prime},\infty)\times Y) $ with $
 r_0 < R^{\prime} < R$, and any $ \epsilon > 0 $ small enough, there exists $
\varphi^{\epsilon} \in C_0^{\infty}({\mathcal M}\setminus
{\mathcal K}) $ such that $ || \varphi - \varphi^{\epsilon} || +
|| A \varphi - A \varphi^{\epsilon} || < \epsilon $. Using the
function $ \theta_{\epsilon} $ defined in Appendix $
\refe{dimension1} $, we set
$$ \varphi^{\epsilon} = \sum_{k } \varphi_k \ast \theta_{\epsilon} \otimes e^{- |\epsilon| \mu_k} \psi_k . $$
It is clearly compactly supported in $ I \times Y $ if $ \epsilon
$ is small enough and smooth since $ \partial_r^j \Delta_h^l
\varphi^{\epsilon} \in L^2 $ for all $j,l \in \Na $. Then, by
Parseval's identity, we  have $ \varphi^{\epsilon} \rightarrow
\varphi $ and using $ \reff{convolution} $ we also have $  A
\varphi^{\epsilon} \rightarrow A \varphi   $. For the last
statement, we first observe that, if $ \varphi $ is compactly
supported, so is $ A \varphi $. We are thus left with the
regularity for which we observe that  $ [\partial_r^j , A_k] =
\sum_{m \leq j} b_{k,m}(r) \partial_r^m $, with $ b_{k,m} $
uniformly bounded by $ \reff{uniff} $, and hence
$$ || \partial_r^j \mu_k^l A_k \varphi_k || \ \leq \ || A_k (\partial_r^j \Delta_h^l \varphi)_k ||
+ C \sum_{m \leq j} || (\partial_r^m \Delta_h^l \varphi)_k || \
\in \ l^2 (\Na_k)
$$
yields the result. \finpreuve

\medskip

Note that the choice of $ C_{B}^{\infty} ({\mathcal M}) $ is
dictated by the following proposition.

\begin{prop} \label{cbinf2} For all $ n \in \Na $, $ z \notin \emph{spec}(H) $
and $ g \in C_0^{\infty}(\Ra) $, we have
$$ \zeta_n (H-z)^{-1} C_B^{\infty} ({\mathcal M}) \subset C_B^{\infty}({\mathcal M})  ,
\qquad \zeta_n g (H) L^2(\mathcal M) \subset
C_B^{\infty}({\mathcal M}) . $$
\end{prop}

\noindent {\it Proof.} This is a direct consequence of standard
elliptic regularity results (see for instance \cite{ChPi,Horm3}),
taking into account the fact that $ \zeta_n = 1 $ near $ \partial
{\mathcal M} $ (if non empty). \finpreuve 

\medskip

We now consider the calculations of $ [H,A] $ and $ [[H,A],A ] $.
Note that these commutators make perfectly sense on $
C_B^{\infty}({\mathcal M}) $ by Propositions $ \refe{cbinf1} $ and
the fact that $ C_B^{\infty}({\mathcal M}) \subset D (H) $.

 We first consider the "free parts", i.e. the commutators involving $ H_0 $
 defined by $ \reff{Hzero}  $.
\begin{prop} \label{warp} There exists $ C $ such that for all  $ R > S > r_0 + 1 $
and all $ \varphi \in C_B^{\infty}({\mathcal M}) $.
\begin{eqnarray}
|| [H_0,A] \varphi || + || [[H_0,A],A] \varphi || \leq C || (H +
i) \varphi ||  . \label{commutator1}
\end{eqnarray}
\end{prop}

\noindent {\it Proof.} Similarly to the proof of Proposition $
\refe{cbinf1} $, we identify $ L^2 ({\mathcal M} \setminus
{\mathcal K}) $ and $ L^2 (I) \otimes L^2 (Y, d \mbox{Vol}_h) $
for notational simplicity. Straightforward calculations show that
\begin{eqnarray}
 \!  \! \! \! \!  i [H_0,A]   \varphi  & = &   \sum_{k \geq 0}  \left(
   2 a_k^{\prime} D_r^2 + 2 a_k \mu_k e^{-2r} - 2 a^{\prime \prime}_{k}
   \partial_r - a_k^{(3)}/2 \right) \varphi_k
   \otimes \psi_k   ,
     \label{freepart1} \\
\! \! \! \! \!  [ [H_0,A],A ]  \varphi  & = & \ \sum_{k \geq 0}
\left( b_k D_r^2 + c_k D_r +d_k \right) \varphi_k
   \otimes \psi_k    \label{freepart2} ,
\end{eqnarray}
where the functions $ b_k(r),c_k(r),d_k(r) $ are given by
\begin{eqnarray}
 b_k = 2 (  a_k a_k^{\prime \prime} - 2 a_k^{\prime \ 2} ) , \qquad c_k = 5 i a_k^{\prime} a_k^{\prime \prime}
 - i a_k a_k^{(3)}  , \qquad  \nonumber \\ d_k = 2 a_k \mu_k e^{-2 r} ( a_k^{\prime} - 2 a_k  )
+ a_k^{\prime} a_k^{(3)} - (a_k a_k^{(4)} - a_k^{\prime \prime \ 2
})/2  .  \nonumber
\end{eqnarray}
One easily checks that $ a_k  \mu_k e^{-2r} $ and $ a_k^2 \mu_k
e^{-2r} $ are uniformly bounded with respect to $ k \in \Na $ and
$ R > S > r_0 + 1 $, thus,  using $ \reff{uniff} $, the result is
direct consequence of the following lemma. \finpreuve

\begin{lemm} \label{FrHiSobolev} For all differential operator $  P $ with coefficients
supported in $ {\mathcal M} \setminus {\mathcal K} $ such that
$$ \tilde{\Psi}^* P \tilde{\Psi}_* = \sum_{j+|\beta| \leq 2}
c_{j,\beta}(r,y) ( e^{-r}D_y )^{\beta} D_r^j  $$ with $
c_{j,\beta} $ bounded on $ I \times U_0 $ for all $ U_0 \Subset U
$ (with the notations of page \pageref{notationcarte}), there
exists $ C $ such that
$$ || P \varphi || \leq C ||(H+i) \varphi ||, \qquad \forall \ \varphi \in D (H) . $$
\end{lemm}

\noindent {\it Proof.} It is a direct application of Lemma 1.3 of
\cite{FrHi1}. \finpreuve

\medskip

We will now give a pseudo-differential approximation of $ A $
which will be useful both for computing the "perturbed parts" $
[A,V] $, $ [A,[A,V]] $ and for the proof of Theorem $ \refe{theo2}
$.

Following \cite{Horm3}, we say that,  for $ m \in \Ra $, $ g \in
S^m (\Ra^{d_1}_x \times \Ra^{d_2}_{\varsigma}) $ if $
|\partial_x^{\alpha}
\partial_{\varsigma}^{\beta} g (x,\varsigma)| \leq C_{\alpha,\beta}
\scal{\varsigma}^{m - |\beta|} $,  for all $ \alpha,\beta $. If $
g \in S^0 (\Ra_r \times \Ra_{\mu}) $ is supported in $ I \times
\Ra $, we clearly define a bounded operator on $ L^2 (I \times Y)
$ by
$$ g (r , \Delta_{h}) \varphi = \sum_{k \geq 0} g (r, \mu_k) \varphi_k \otimes \psi_k . $$
Abusing the notation for convenience, we still denote by $ g (r,
\Delta_{h}) $ the pullback of this operator on $ L^2 ({\mathcal
M}\setminus {\mathcal K} ) $, extended by $ 0 $ on $ L^2
({\mathcal K}) $. If $ \theta \in C^{\infty} (Y) $, we also denote
by $ \theta $ (instead of $ 1 \otimes \theta $) its natural
extension to $ I \times Y $ which is independent of $ r $. Our
pseudo-differential approximation of $A$ will mainly follow from
the following result.

\begin{prop} \label{approxpseudo} Let $ g \in S^{0}(\Ra_r \times \Ra_{\mu}) $ be supported in $  I \times Y $.
For all coordinate patch $ U_Y \subset Y $,  all $ \theta,
\tilde{\theta} \in C^{\infty}_0 (U_Y) $  such that $
\tilde{\theta} \equiv 1 $ near the support of $ \theta $ and all $
N $ large enough, there exists $ g_N \in S^{0} ( \Ra_{r,y}^{n}
\times \Ra^{n-1}_{\eta}) $ and an operator $ {\mathcal
R}_N^{\theta} : L^2 (\mathcal M) \rightarrow L^2 (\mathcal M) $
such that
\begin{eqnarray}
 \theta g (r,\Delta_{h}) = G_N^{\theta} + {\mathcal R}_N^{\theta}
 \label{compright}
\end{eqnarray}
 where $ G_N^{\theta} = \tilde{\Psi}_* \left( (\Psi^* \theta) (y)
g_N (r,y,D_y) (\Psi^* \tilde{\theta})  (y) \right) \tilde{\Psi}^*
$ (with the notation $ \reff{chart} $) and
\begin{eqnarray}
 \left| \left|
\Delta_{h}^j {\mathcal R}_{N}^{\theta} \Delta_{h}^k \varphi \right|
\right| \leq C_{j,k} || \varphi ||, \qquad \varphi \in
C_c^{\infty}({\mathcal M}) \label{regularisation} , \\
\left| \left| \Delta_{h}^j [ D_r , {\mathcal R}_{N}^{\theta} ]
\Delta_{h}^k \varphi \right| \right| \leq C_{j,k} || \varphi ||,
\qquad \varphi \in C_c^{\infty}({\mathcal M}), \label{comDr}
\end{eqnarray}
for all $ j ,k \leq N $. If $ p_h $ is the principal symbol of $
\Delta_{h} $, we actually  have
$$ g_{N}(r,y,\eta) = g (r,p_h (y,\eta)) + \sum_{ 1 \leq
j \leq j_N  } \sum_{l} d_{jl}(y,\eta) \partial_{\mu}^j
g(r,p_h(y,\eta))
$$
where $ d_{jl} $ are polynomials of degree $ 2j-l $ in $ \eta $,
obtained as universal sums of products of the full symbol of $
\Delta_h $ in coordinates $ (y,\eta) $. 

More generally, if $
(g_{\lambda})_{\lambda \in \Lambda} $ is a  bounded family in $ S^0
(\Ra_r \times \Ra_{\mu}) $ with support in $ I \times \Ra $, the
associated family $ (g_{\lambda,N})_{\lambda \in \Lambda} $ is
bounded in $ S^0 $ and the constant $ C_{j,k} $ in $
\reff{regularisation} $ can be chosen independent of $ \lambda \in
\Lambda $.
\end{prop}

The proof  is given in Appendix $ \refe{functional} $. Note that,
strictly speaking, this proposition is not a direct consequence of
 the standard functional
calculus for elliptic pseudo-differential operators on closed
manifolds \cite{Seel1} since $g$ depends on the extra variable $ r
$. However, the proof follows from minor adaptations of the
techniques of \cite{HeRo1,Seel1}.

\smallskip

\noindent {\bf Remark 1.} The operators $ g (r,\Delta_{h}) $ and $
G_N^{\theta} $ commute with operators of multiplication by
functions of $ r $, hence so does $ {\mathcal R}_N^{\theta} $.
\newline
{\bf Remark 2.} In $ \reff{regularisation} $, we have abused the
notation by identifying $ \Delta_{h}  $, which acts on functions on
$ Y$, with its natural extension acting on functions on $
{\mathcal M} $ which are supported in $ {\mathcal M} \setminus {\mathcal K} $.

\medskip

The previous proposition is motivated by the fact that we can
write
\begin{eqnarray}
 A = g_{R,S}(r,\Delta_{h}) r D_r + \tilde{g}_{R,S} (r,\Delta_{h}) .
 \label{formedeA}
\end{eqnarray}
with functions  $ g_{R,S} $ and $ \tilde{g}_{R,S} $ belonging to $
S^0 (\Ra_r \times \Ra_{\mu}) $ as explained  by the following lemma.
\begin{lemm} There exist two families $ g_{R,S} ,  h_{R,S}  \in  
S^0 (\Ra_r \times \Ra_{\mu}) $, bounded for $ R > S > r_0 +1 $, 
supported in $ r > R $ and such that
$$ g_{R,S} (r,\mu_k) = a_k (r) /r , \qquad \tilde{g}_{R,S} (r,\mu_k) =  a^{\prime}_k (r) / 2 i $$
for all $ k \geq 0 $.
\end{lemm}

\noindent {\it Proof.} With $ \gamma \in C^{\infty}(\Ra_{\mu}) $
such that $ \gamma = 1 $ on  $ \Ra^+ $ and $ \mbox{supp} \ \gamma
\in [-1/2,\infty ) $, we may choose
$$ g_{R,S}(r,\mu) = \gamma (\mu) \chi_{R}(r) \left( 1 + 2 \frac{S}{r} - \frac{1}{2r} \log (1+\mu)
 \right) \xi_S \left( r - \frac{1}{2} \log (1 + \mu) \right) . $$
It is easily seen to belong to $ S^0 (\Ra^2) $ and the boundedness
with respect to $ R,S$ follows from
$$ \partial_{\mu}^j \xi_S \left( r- \frac{1}{2}\log (1+\mu) \right) = \sum_{1 \leq k \leq j} c_{jk} S^{-k}
 \xi^{(k)} \left( \frac{r}{S} -
\frac{1}{2S} \log(1+\mu) \right) (1 + \mu)^{-j}, $$  the fact that
$  - S /2 \leq r - \frac{1}{2} \log (1+\mu) \leq r + \log 2^{1/2}
$ on the support of $ \gamma(\mu)  \xi_S (r -
\frac{1}{2}\log(1+\mu)) $ and the fact that $ S /r $ is bounded on
the support of $ \chi_R (r) $. Then, we may choose $
\tilde{g}_{R,S} = g_{R,S} + r \partial_r g_{R,S} $ since one
checks similarly that $ r \partial_r g_{R,S} $ is bounded in $ S^0
$. \finpreuve

\medskip

We are now ready to study the contribution of the perturbation $ V
$ for the commutators.
 \begin{prop} \label{souscomm} There exists $ C > 0 $
such that, for all  $ R > S > r_0  $ and  all $ \varphi \in
C_B^{\infty}({\mathcal M}) $
\begin{eqnarray}
 || \scal{r}^j [A,V] \varphi || & \leq &  C
R^{j-1}||(H+i)\varphi||, \qquad j = 0 , 1 , \label{souscomm1} \\
 ||  [A,[A,V]] \varphi || & \leq & C  ||(H+i)\varphi|| . \label{souscomm2}
\end{eqnarray}
\end{prop}

\noindent {\it Proof.} Dropping the subscripts $ R,S $ on $ g $
and $ \tilde{g} $, we have
$$ A = g (r,\Delta_{h}) r D_r + \tilde{g}(r,\Delta_{h}) = \left( G_N + {\mathcal R}_N \right) r D_r
+ \tilde{G}_N + \tilde{\mathcal R}_N $$ with $ G_N = \sum_l
G_N^{\theta_l} $ and $ {\mathcal R}_N =\sum_l {\mathcal
R}_N^{\theta_l} $ associated to $ g $ by mean of proposition $
\refe{approxpseudo} $ and of a partition of unit  $ \sum_{l}
\theta_{l} = 1 $ on $ Y$. Of course, $ \tilde{G}_N $ and $
\tilde{\mathcal R}_N $ are similarly associated to  $ \tilde{g} $.
Note that $ g (r,\Delta_{h}) $ and $ G_N $ map $
C_c^{\infty}(\mathcal M) $ into $ C_0^{\infty}({\mathcal M}
\setminus {\mathcal K}) $ and thus so does $ {\mathcal R}_N $.
Therefore, on $ C_B^{\infty} ({\mathcal M}) $, we have
$$ [A,V] =  \left( [G_N, V ] + [{\mathcal R}_N,V ] \right) r D_r
+  g (r,\Delta_{h}) [r D_r,V] + [\tilde{G}_N,V] + [\tilde{\mathcal
R}_N , V ] . $$ We study the terms one by one. Note first  that
$ [G_N,V]rD_r = r \scal{r}^{-2} [G_N , \scal{r}^2 V] D_r$. If $
\tilde{\Psi}_l $ is  associated to a coordinate chart $ \Psi_l $
defined in a neighborhood of $ \mbox{supp} \ \theta_l $ by $
\reff{chart} $, we have $$  \tilde{\Psi}_l^* [G_N^{\theta_l} ,
\scal{r}^2 V] \tilde{\Psi}_{l  *} = \sum_{|\beta| \leq 1}
q_{\beta}(r,y,D_y) (e^{-r}D_y)^{\beta}
 $$
 with $ q_{\beta} \in S^0 $ which depends, in a bounded way, on $ R > S > r_0 $
 and is supported in $ r \geq R $. This follows by standard pseudo-differential
 calculus and thus, by Lemma $ \refe{FrHiSobolev} $, we have
$$ \left| \left| \scal{r}^j [G_N,V]rD_r \varphi \right| \right| \leq C
R^{j-1}|| (H+i)  \varphi
||, \qquad \varphi \in C_B^{\infty}({\mathcal M})
 $$
with $C $ independent of $ R >S > r_0  $. Similarly, we get the
same estimate for $ [{\mathcal R}_N,V] r D_r $ since $ \scal{r}^2
[{\mathcal R}_N,V] $ is a bounded operator, uniformly with respect
to $ R > S > r_0  $, with range supported in $ r \geq R $. The
same holds for $ [\tilde{G}_N,V]  $ and $  [\tilde{\mathcal R}_N ,
V ] $. Finally, $ \scal{r} [V,rD_r] $ is an operator of the form
considered in Lemma $ \refe{FrHiSobolev} $, whereas $ ||
\scal{r}^{-1} g (r,\Delta_{h}) || \leq C R^{-1} $, so $
\reff{souscomm1} $ follows.

We now consider $ [A,[A,V]] $. We only study $
[g(r,\Delta_{h})rD_r,[g(r,\Delta_y)rD_r,V]] $, since the other terms
can be studied similarly and involve less powers of $ r D_r $.
This double commutator reads
\begin{eqnarray}
 \left[ g (r,\Delta_{h}) , [g(r,\Delta_y)rD_r,V] \right] r D_r + g
(r,\Delta_{h}) \left[ r D_r , [g(r,\Delta_y)rD_r,V] \right] = \qquad
\qquad \nonumber \\ \qquad \left[ G_N , [G_N rD_r,V] \right] r D_r
+ G_N \left[ r D_r , [G_N rD_r,V] \right] + I_N D_r^2 + J_N D_r +
K_N
\end{eqnarray}
where $ I_N , J_N , K_N $ are bounded operator on $ L^2 ({\mathcal
M }) $, uniformly with respect to $ R > S > r_0 + 1 $. This
clearly follows from Proposition $ \refe{approxpseudo} $ and the
fact that $ 1 \otimes (\Delta_{h}+1)^{-j} ( r^2 V )   1 \otimes
(\Delta_{h} + 1) ^{-k}  $ is bounded if $ j + k \geq 1 $. Precisely,
 $ 1 \otimes (\Delta_{h} +1)^{-1}$ is actually defined on $ L^2 (I
\otimes Y) $ but, here, it is identified with its pullback on $
L^2 ({\mathcal M}\setminus {\mathcal K}) $. By Lemma $
\refe{FrHiSobolev} $, $ ||(I_N D_r^2 + J_N D_r + K_N ) \varphi ||
\leq C ||(H+i) \varphi || $. On the other hand, for all $
\theta_{l_1} $ and $ \theta_{l_2} $ associated with overlapping
coordinate patches, we have
$$  \tilde{\Psi}_{l_1}^* \left[ G_N^{\theta_{l_1}} , [G_N^{\theta_{l_2}} rD_r,V] \right]
r D_r \tilde{\Psi}_{l_1 *} =   \sum_{|\beta| + k \leq 2}
\tilde{q}_{\beta}(r,y,D_y) (e^{-r}D_y)^{\beta} D_r^k $$ with $
\tilde{q}_{\beta} $ bounded in $ S^0 $ for $ R > S > r_0  $. This
follows again from the usual composition rules of
pseudo-differential operators and it clearly implies that
$$ \left| \left| \left[ G_N , [G_N rD_r,V] \right]
r D_r \varphi \right| \right| \leq C ||(H+i)\varphi ||, \qquad
\varphi \in C_B^{\infty}({\mathcal M}),
$$
with $C $ independent of $ R , S $. Similarly, the same holds for
$ G_N \left[ r D_r , [G_N rD_r,V] \right] $ and the result
follows. \finpreuve

\medskip

We conclude this subsection with the following proposition which
summarizes what we know so far on $A$ and $H$.

\begin{prop} \label{checkconditions} With $ {\mathcal D} = C_B^{\infty}({\mathcal M}) $,
all the conditions from $ \reff{coreA} $ to $ \reff{comm2} $ hold.
Furthermore, in $ \reff{comm2} $, $ C_{H,A} $  can be chosen
independently of $ R > S
> r_0 + 1 $.
\end{prop}

\noindent {\it Proof.} Using Lemma $ \refe{FrHiSobolev} $, it is
clear that $ [ H , \zeta_n ] \rightarrow 0 $ strongly on $ D (H) $
as $ n \rightarrow \infty $. Therefore, all the conditions from $
\reff{coreA} $ to $ \reff{strongH} $ are fulfilled. In particular,
$ C_B^{\infty} ({\mathcal M}) $ is a core for $ H $, hence
Propositions $ \refe{warp} $ and $ \refe{souscomm} $ yield the
existence of $ [H,A]^0 $, and thus  $ \reff{comm1} $ and $
\reff{comm2} $ hold. It only remains to prove $ \reff{stable} $.
Assume for a while that, for all $ \varphi , \psi \in C_B^{\infty}
({\mathcal M}) $,
\begin{eqnarray}
 ((H-z)^{-1} \varphi , A \psi )  - (A \varphi , (H-\bar{z})^{-1} \psi)  = \left(
(H-z)^{-1}[H,A]^0 (H-z)^{-1} \varphi, \psi \right) .
\label{commform}
\end{eqnarray}
Then this holds for all $ \varphi , \psi \in D (A) $. Since $
(H-z)^{-1}[H,A]^0 (H-z)^{-1} $ is bounded, $ \reff{commform} $
yields
$$ \left| ((H-z)^{-1} \varphi , A \psi) \right| \leq C (|| A \varphi || + || \varphi ||)||\psi|| , $$
which shows that $ (H-z)^{-1} \varphi \in D (A^*) = D (A) $ for
all $ \varphi \in D (A) $ and hence $ \reff{stable} $. Let us show
$ \reff{commform} $. By $ \reff{regularity}  $, the right hand
side of $ \reff{commform} $ can be written as the limit, as $ n
\rightarrow \infty $, of $ ( [H,A] \zeta_n (H-z)^{-1} \varphi ,
\zeta_n (H-\bar{z})^{-1} \psi ) $ i.e. the limit of
\begin{eqnarray}
 \left(  \zeta_n (H-z)^{-1} \varphi , A \zeta_n \psi + A [H,
\zeta_n ] (H- \bar{z})^{-1} \psi \right)
 -  \left( A \rho_n \varphi + A [H, \zeta_n ] (H-z)^{-1}
\varphi  , \zeta_n (H-\bar{z})^{-1} \psi \right)  . \nonumber
\end{eqnarray}
By  $ \reff{formedeA} $, Lemma $ \refe{FrHiSobolev} $ and the fact
that $ 2^{-n } r \zeta^{\prime} (2^{-n} r) \rightarrow 0$, it is
clear that $ A [H,\zeta_n] (H-z)^{-1} \varphi \rightarrow 0 $. The
same holds for $ \psi $ of course and thus $ ( [H,A] \zeta_n
(H-z)^{-1} \varphi , \zeta_n (H-\bar{z})^{-1} \psi ) $ converges
to the left hand side of $ \reff{commform} $. This completes the
proof. \finpreuve
\subsection{Positive commutator estimate}

This subsection is devoted to the proof of a positive commutator
estimate of the form $ \reff{commpositif} $ at large energies $
\lambda $ (with control on $ \delta $ with respect to $ \lambda
$).

We start with some notation. Let $ \Xi_{R,S} $ be the pullback on
$ L^2 ({\mathcal M} \setminus {\mathcal K}) $ (extended by $0$ on
$ L^2({\mathcal K}) $) of the operator defined on $ L^2 (I \times
Y ) $ by
$$ \varphi \mapsto \sum_{k \geq 0}  \chi_R^{1/2} (r) (1-\xi_S^{1/2})(r-\log \nu_k) \varphi_k \otimes \psi_k $$
with the notation $ \reff{coefFourier} $. We also set $
\widetilde{\Xi}_{R,S} = \chi_R^{1/2} - \Xi_{R,S} $. 
 Similarly, $ \Xi_{R,S}^{\prime}, \Xi_{R,S}^{\prime \prime} $ are the
operators respectively  defined by $ \partial_r \left( \chi_R^{1/2} (r)
(1-\xi_S^{1/2})(r-\log \nu_k) \right)$ and $ \partial_r^2 \left(
\chi_R^{1/2} (r) (1-\xi_S^{1/2})(r-\log \nu_k) \right) $.

\begin{prop} \label{commutatorestimate}
 There exists $ C $ such that, for all $ \lambda \gg 1 $, all $ F \in
 C_0^{\infty}([0,2\lambda], [0,1]) $ and all $ R > S > r_0 + 1 $, one has
\begin{eqnarray}
 F (H) i [H,A]^0 F (H) - 2 H F (H)^2 \qquad \geq \qquad \qquad \qquad \qquad \qquad \qquad \qquad
 \qquad \qquad \qquad
  \nonumber \\
\qquad \qquad \qquad  -  C \lambda \left( || F (H) \scal{r}^{-1}||
+ || F (H) (\chi_R - 1)|| + || F (H) (1 - \widetilde{\Xi}_{R,S}^2
) || + S^{-1} + \lambda^{-1} \right) . \label{shrinkfonction}
\end{eqnarray}
\end{prop}

\noindent {\it Proof.}
 We first note that the right hand side of $
\reff{freepart1} $ is nothing but $ 2 D_r a_{k}^{\prime} D_r + 2
a_k \mu_k e^{-2 r} - a_k^{(3)} / 2 $. Since $ a_k^{\prime}(r) \geq
\chi_R (r) \xi_S (r - \log \nu_k) $ and  $ a_k (r) \geq S \chi_R
(r) \xi_S (r - \log \nu_k) $, we get
$$ i [H_0,A] \geq 2
D_R \widetilde{\Xi}_{R,S}^2 D_R + 2 \widetilde{\Xi}_{R,S} e^{-2r}
\Delta_h \widetilde{\Xi}_{R,S} - C S^{-2} . $$
This estimates, as well as the following, holds 
when tested against elements of $ {\mathcal D} = C_B^{\infty}({\mathcal M}) $.
For any $
a \in C^{\infty}(\Ra) $, one has $ D_r a^2 D_r = a D_r^2 a + a
a^{\prime \prime} $, so the above inequality yields
$$ i [H_0,A] \geq 2 \widetilde{\Xi}_{R,S} H_0 \widetilde{\Xi}_{R,S} -
(n-1)^2 \widetilde{\Xi}^2_{R,S} /4 -  C S^{-2} , $$ for $ R
> S > r_0 +1 $. We then write
\begin{eqnarray}
\widetilde{\Xi}_{R,S} H_0 \widetilde{\Xi}_{R,S} =H_0 + ( \chi_R -
1 ) H_0 + (1-\widetilde{\Xi}_{R,S}^2)H_0 -
\widetilde{\Xi}^{\prime}_{R,S} \widetilde{\Xi}_{R,S}
\partial_r - \widetilde{\Xi}_{R,S} \widetilde{\Xi}^{\prime
\prime}_{R,S}\nonumber
\end{eqnarray}
and this implies that, on $ D (H) $,  $ i [H,A]^0 \geq 2 H +
Q_{R,S} - C  $ with $ C$ independent of $ R,S $ and  $$ Q_{R,S} =
i [V,A]^0 - V + ( \chi_R - 1 ) H_0 + (1-\widetilde{\Xi}_{R,S}^2)H_0
- \widetilde{\Xi}^{\prime}_{R,S} \widetilde{\Xi}_{R,S}
\partial_r - \widetilde{\Xi}_{R,S} \widetilde{\Xi}^{\prime
\prime}_{R,S},
$$
where $ [V,A]^0 $ is the closure of $ [V,A] $ (defined on
$C_B^{\infty}({\mathcal M})$) on $ D (H) $ . Then, using 
Lemma  $ \refe{FrHiSobolev} $, we have  $ || H_0 F (H) || + ||
\chi_R^{1/2} \partial_r F (H) ||+ || \scal{r}^{2} V F (H) || \leq
C \lambda $, and using Proposition $ \refe{souscomm} $, the result
follows. \finpreuve

\medskip

Note that, if $ F $ is supported close enough to $ \lambda $, $ 2
H F (H)^2 \geq 3 \lambda F^2 (H) /2 $ and thus we will get $
\reff{commpositif} $ by making the bracket of the right hand side
of $ \reff{shrinkfonction} $ small enough.

\medskip

 Using the technique of
\cite{FrHi1}, we are able to estimate $ || F (H) (1 -
\widetilde{\Xi}_{R,S}^2 ) || $ for suitable $F$. Let us recall
the proof of this fact. For $ R > S > r_0 + 1 $, a direct calculation yields
$$ \Xi_{R,S}^2 H_0 + H_0 \Xi_{R,S}^2 = 2 \Xi_{R,S} H_0 \Xi_{R,S} - 2 \left( \Xi_{R,S}^{\prime} \right)^2 .  $$
On the other hand, $ e^{-2r} \mu_k \chi_R (r) \geq e^{S} - e^{-2R}
$ on the support of $ \chi_R (r) \xi_S (r - \log \nu_k) $ so we
also have $ \Xi_{R,S} H_0 \Xi_{R,S} \geq \left( e^{S} - e^{-2R}
\right) \Xi_{R,S}^2 $, and we obtain
$$ \Xi_{R,S}^2 \left( \tau(H_0 - \lambda) - z \right) + \left( \tau ( H_0 - \lambda ) - \bar{z} \right)
 \Xi_{R,S}^2 \geq 2 \tau  \left(e^{S} - e^{-2R} - \lambda - \frac{\mbox{Re}z}{\tau} \right)
 \Xi_{R,S}^2 - 2 \tau \left( \Xi_{R,S}^{\prime} \right)^2 , $$
 for all real $ \tau \ne 0 $, $ z \in \Ca $ and $ \lambda \in \Ra
 $.
 Testing this inequality against $ \left( \tau (H_0 - \lambda) - z \right)^{-1} \psi
 $, we  get
\begin{eqnarray}
  2 \left( \psi , \Xi_{R,S}^2 \left( \tau (H_0 - \lambda) - z \right)^{-1} \psi  \right)
 + 2 \tau \left| \left| \Xi_{R,S}^{\prime} \left( \tau (H_0 - \lambda)-z \right)^{-1} \psi
  \right| \right|^2  \qquad \qquad \nonumber \\
  \qquad \qquad \geq 2 \tau \left(e^{S} - e^{-2R} - \lambda - \frac{\mbox{Re}z}{\tau} \right)
 \left| \left| \Xi_{R,S} \left( \tau (H_0 - \lambda)-z \right)^{-1} \psi
  \right| \right|^2 \nonumber
  \end{eqnarray}
  and this clearly implies, provided $ e^{S} - e^{-2R} - \lambda - \mbox{Re}z/\tau > 0
  $, $ \tau > 0 $ and $ R > S > r_0 + 1 $, that
\begin{eqnarray}
  \left| \left| \Xi_{R,S} \left( \tau (H_0 - \lambda)-z
\right)^{-1}
  \right| \right| \leq  \frac{1}{|\mbox{Im}z|}
  \left( e^{S} - e^{-2R} - \lambda - \frac{\mbox{Re}z}{\tau}
  \right)^{-1/2}\! \!
  \left( \frac{ C_{\chi,\xi} }{S}+
  \frac{|\mbox{Im}z|^{1/2}}{\tau^{1/2}} \right) . \label{semiclassique}
\end{eqnarray}
This estimate is essentially taken from \cite{FrHi1} and is the
main tool of the proof of
\begin{prop} \label{HelfferSjostrand} Let $ F_0 \in C_0^{\infty} ([-1,1],\Ra) $ such that $ 0 \leq F_0 \leq 1
$. There exists $ C $ such that, with $$ R = \log 5 \lambda ,
\qquad S = \log 4 \lambda , \qquad \tau = \lambda^{-1} , $$ we
have
$$ \left| \left| F_0 \left( \lambda^{-1} H - 1 \right) (1 - \widetilde{\Xi}_{R,S}^2) \right| \right|
\leq C \lambda^{-1/2} ( \log \lambda )^{ -1 }   , \qquad \lambda
\gg 1 .
$$

\end{prop}

\noindent {\it Proof.}  We shall use Helffer-Sj\"ostrand formula (see for instance \cite{DiSj1}) , i.e. 
$$ F_0 \left( \tau (H_0 - \lambda) \right) = \frac{1}{2  \pi} \int \!
\! \! \int_{\Ra^2} \overline{\partial} \tilde{F_0} (u+iv) \left(
\tau (H_0-\lambda) - u - i v \right)^{-1}
 \mbox{d}u  \mbox{d} v , $$
 where $ \overline{\partial} = \partial_u +
i \partial_v$, 
$ \tilde{F_0} \in C_0^{\infty} (\Ca) $ is such that $ \tilde{F_0}_{| \Ra} = F_0
 $ and $ \overline{\partial} \tilde{F_0}   = {\mathcal O}(|v|^{\infty})
 $ near $ v = 0 $. As a direct consequence of $ \reff{semiclassique}
 $ with $ R = \log 5 \lambda $, $ S = \log 4 \lambda $, $ \tau = \lambda^{-1} $, and assuming
 that $ |\mbox{Re}z | \leq 2 $ on the support of $ \tilde{F_0} $,
 Helffer-Sj\"ostrand formula gives
$$  \left| \left| F_0 \left( \lambda^{-1} H_0 - 1 \right) (1 - \widetilde{\Xi}_{R,S}^2) \right| \right|
\leq C_{F_0} \lambda^{-1/2} (\log \lambda)^{-1}, \qquad \lambda
\gg 1 .
$$ We are thus left with the study of $ ||   \left( F_0 \left(
\lambda^{-1} H - 1 \right) - F_0 \left( \lambda^{-1} H_0 - 1
\right) \right) (1 - \widetilde{\Xi}_{R,S}^2)  || $ or,
equivalently, with $ ||   (1 - \widetilde{\Xi}_{R,S}^2) \left( F_0
\left( \lambda^{-1} H - 1 \right) - F_0 \left( \lambda^{-1} H_0 -
1 \right) \right) || $. Using the resolvent identity
$$ \left(
\lambda^{-1} H -1 - z \right)^{-1} - \left( \lambda^{-1} H_0- 1 -
z \right)^{-1} = \left( \lambda^{-1} H_0 - 1 - z \right)^{-1}
\lambda^{-1} V \left(  \lambda^{-1} H- 1 - z \right)^{-1} ,
$$ and the fact that $ V $ is $H $ bounded with relative bound $ < 1 $,
which implies that, for some $ C $ independent of $ \lambda \gg 1
$ and $ z \in \mbox{supp} \ \tilde{F_0} $,  $ || \lambda^{-1} V
\left( \lambda^{-1} H- 1 - z \right)^{-1} || \leq C | \mbox{Im} z
|^{-1} $, another application of Helffer-Sjostrand formula implies
that
$$ \left| \left|   (1 - \widetilde{\Xi}_{R,S}^2) \left( F_0 \left( \lambda^{-1}
H - 1 \right) - F_0 \left( \lambda^{-1} H_0 - 1 \right) \right)
\right| \right| \leq  C_{F_0} \lambda^{-1/2} (\log \lambda)^{-1},
\qquad \lambda \gg 1 .
$$
The result follows. \finpreuve

\medskip

We can now explain how to get an estimate of the form $
\reff{commpositif} $. For any fixed $ 0 < \epsilon < 1 $,  one can
clearly choose $ F_0 \in C_0^{\infty}(\Ra,\Ra) $ as above such
that for all $ F_1 \in C_0^{\infty} (\Ra,\Ra) $ supported in $ [(1
- \epsilon) \lambda  , (1 + \epsilon) \lambda ] $, we have $ F_1
(E) = F_0 (\lambda^{-1} E - 1) F_1 (E) $ for all $ E \in \Ra $.
Thus, for all such $ F_1 $ satisfying $  0 \leq F_1 \leq 1 $,
Propositions $ \refe{commutatorestimate} $ and $
\refe{HelfferSjostrand} $ imply that, for $ \lambda \gg 1 $,
\begin{eqnarray}
 F_1 (H) i [H,A]^0 F_1 (H) \geq (2 - 2 \epsilon) \lambda F_1 (H)^2 - \qquad \qquad \qquad \qquad
\qquad \qquad \qquad \qquad  \nonumber \\
 \qquad \qquad \qquad \qquad \qquad \qquad C
\lambda \left( || F_1 (H) \scal{r}^{-1}|| + || F_1 (H)
(\chi_R-1)|| + ( \log \lambda )^{-1} \right),  \label{commpresque}
\end{eqnarray}
 if $ R = \log  5
\lambda  $ and $ S = \log 4 \lambda $.
 Then, if we assume  that there exists $ 0< s_0 \leq 1  $ such that
\begin{eqnarray}
 || \scal{r}^{-s_0} (H - \lambda \pm i 0) \scal{r}^{-s_0} || \leq
\varrho (\lambda) , \qquad \lambda \gg 1 , \label{Vodev}
\end{eqnarray}
with $ \varrho (\lambda) > \lambda^{-1}/C $, we can choose $ F_1 $
in view of  Lemma $ \refe{easytrick} $. Indeed,  $ \chi_R - 1 $ is
supported in $ |r| \leq C \log \lambda $, so we have $ || F_1 (H)
(\chi_R-1)|| \leq C || F_1 (H) \scal{r}^{-s_0}|| (\log
\lambda)^{s_0}  $, and thus $ \reff{commpresque} $ reads
$$   F_1 (H) i [H,A]^0 F_1 (H) \geq \frac{3}{2} \lambda F_1 (H)^2 - C
\lambda \left( (\log \lambda)^{s_0} || F_1 (H) \scal{r}^{-s_0}|| +
( \log \lambda )^{-1} \right) . $$ Hence,  if $ F_1  $ supported
in $ [ \lambda - c \varrho (\lambda)^{-1} (\log \lambda)^{-2 s_0}
, \lambda + c \varrho (\lambda)^{-1} (\log \lambda)^{-2 s_0} ] $
with $ c
> 0 $ small enough (independent of $ \lambda $), Lemma $
\refe{easytrick} $ clearly shows that
$$  F_1 (H) i [H,A]^0 F_1 (H) \geq (2 - 2 \epsilon) \lambda F_1 (H)^2 - \lambda / 2, \qquad \lambda \gg 1 . $$
Note that the condition  $ [ \lambda - c \varrho (\lambda)^{-1}
(\log \lambda)^{-2 s_0} , \lambda + c \varrho (\lambda)^{-1} (\log
\lambda)^{-2 s_0} ] \subset [(1 - \epsilon) \lambda  ,
(1+\epsilon) \lambda ] $ is ensured, for $ \lambda \gg 1 $ , by
the fact that $ \varrho (\lambda) \geq \lambda^{-1}/C $. All this
easily leads to the
\begin{theo} \label{checkconditions2} Let $ A_{\lambda} $ be the operator given in
Definition $ \refe{definition} $, with $ R = \log (5 \lambda) $
and $ S = \log (4 \lambda) $. Assume that $ \reff{Vodev} $ holds
for some $ 0 < s_0 \leq 1 $ and $ \varrho (\lambda) > \lambda^{-1}
/C $ and let
$$ f_{\lambda} (E) =  f \left( \frac{E - \lambda}{\delta_{\lambda}} \right) , \qquad
\delta_{\lambda}= (\log \lambda)^{- 2 s_0} \varrho (\lambda)^{-1}
/ C , $$
 with $ f \in C_0^{\infty}(\Ra,[0,1]) $, supported in $ [-3,3] $
 and  $ f = 1 $ on $ [-2,2] $. Then, for $ C $  large
 enough, we have
 $$   f_{\lambda} (H) i [H,A_{\lambda}]^0 f_{\lambda} (H) \geq  \lambda f_{\lambda} (H)^2 ,
  \qquad \lambda \gg 1 .  $$
\end{theo}

\section{Proofs of the main results} \label{application}
\setcounter{equation}{0}
\subsection{Proof of Theorem $ \bf \refe{theo1} $}
By Proposition $ \refe{checkconditions} $ and Theorem $
\refe{checkconditions2} $, we are in position to use Theorem $
\refe{theoclef} $. Here the parameter $ \nu $ is $ \lambda $ and we consider
$$ H_{\nu} = H, \qquad A_{\nu} = A_{\lambda}/\lambda^{1/2}, \qquad \alpha_{\nu}
=\lambda^{1/2} , \qquad \delta_{\nu} = (\log \lambda)^{- 2 s_0}
\varrho (\lambda)^{-1} / C , $$ with $ C $ large enough,
independent of $ \lambda $. Assuming that $ \varrho (\lambda) \geq
\lambda^{-1/2} / C $ ensures that $  \delta_{\nu}
\alpha_{\nu}^{-1} \leq 1 $. Using the forms of $ C_0 , C_{1/2},
C_1 $ given on page \pageref{constantespages}, it is easy to check
that
$$ C_{0,\nu} \leq C \alpha_{\nu} \delta_{\nu}^{-2} , \qquad C_{1/2,\nu} \leq C \alpha_{\nu}^{1/2}
\delta_{\nu}^{-1} , \qquad C_{1,\nu} \leq C \alpha_{\nu}
\delta_{\nu}^{-1} .
$$ Furthermore, it is clear that, with $ f_{\nu} = f_{\lambda} $,
we have $ ||[H_{\nu},A_{\nu}] f_{\nu} (H_{\nu}) || \leq C
\alpha_{\nu} $, so Theorem $ \refe{theoclef} $ yields
$$ \left| \left| \scal{A_{\lambda}/\lambda^{1/2}}^{-s} (H - \lambda \pm i 0)^{-1}
  \scal{A_{\lambda}/\lambda^{1/2}}^{-s} \right| \right| \leq   C \varrho
  (\lambda)^{-1} (\log \lambda)^{ 2 s_0} , \qquad \lambda \gg 1 . $$
  Then, by writing
$$ (H - z)^{-1} =(H-Z)^{-1} +
 (z-Z)(H-Z)^{-2} + (z-Z)^2 (H-Z)^{-1}(H-z)^{-1}(H-Z)^{-1} $$
with $ Z = \lambda + i \lambda^{1/2} $, $ z = \lambda \pm i
\varepsilon $  and letting $  \varepsilon \rightarrow 0$, Theorem
$ \refe{theo1} $ will be a consequence of the following
lemma.
\begin{lemm} There exists $ C_s > 0 $ such that
$$ \left| \left| W_{-s} (H- \lambda - i \lambda^{1/2})^{-1} \scal{A_{\lambda}/\lambda^{1/2}}^s \varphi
\right| \right| \leq C_s \lambda^{-1/2} (\log \lambda)^{s} ||
\varphi ||, \qquad \varphi \in D (A_{\lambda}) , \ \ \lambda \gg 1
.
 $$
\end{lemm}

\noindent {\it Proof.} We follow \cite{PSS1}, i.e. argue by
complex interpolation. We only have to consider the case $ s = 1 $
and thus study $ \lambda^{-1/2} W_{-1} (H-\lambda - i
\lambda^{1/2})^{-1} A_{\lambda}$ which we can write, on $ D
(A_{\lambda}) $, as $$ \lambda^{-1/2} W_{-1} A_{\lambda} (H-\lambda -
i \lambda^{1/2})^{-1} -  \lambda^{-1/2} W_{-1} (H-\lambda - i
\lambda^{1/2})^{-1} [ H , A_{\lambda}]^0 (H-\lambda - i
\lambda^{1/2})^{-1} . $$ The second term is $ {\mathcal
O}(\lambda^{-1/2}) $ since $ [H,A_{\lambda}]^0 (H+i)^{-1} $ is
uniformly bounded by Propositions $ \refe{warp} $ and $
\refe{souscomm} $, and $ || (H+i)(H -  \lambda - i
\lambda^{1/2})^{-1} || = {\mathcal O} (\lambda^{-1/2}) $. For the
first term, it is easy to check that $ || \chi_{r_0 + 1} D_r (H -
\lambda - i \lambda^{1/2})^{-1} || \leq C $, using Proposition $
\refe{FrHiSobolev} $ and thus
$$ \left| \left| \lambda^{-1/2} W_{-1} A_{\lambda} (H-\lambda
- i \lambda^{1/2})^{-1} \right| \right| \leq C \lambda^{-1/2}
\sup_{k \geq 0,  \ r \geq R}\left( 1 + \frac{ (r + 2 S - \log
\nu_k) \chi_{R}(r)\xi_{S}(r-\log \nu_k)}{w (r - \log \nu_k) }
\right)
 $$
 with $ R = \log (5 \lambda) $ and $ S = \log (4 \lambda ) $. It
 is not hard to check that the supremum is dominated by $ C \log \lambda
 $ and the result follows. \finpreuve

\subsection{Proof of Theorem $ \bf \refe{theo2} $}

We first prove that $  w (r - \log \scal{\eta}) $ is a temperate
weight, i.e. satisfies $ \reff{temperateweight} $ below.

\begin{lemm} \label{poids}
 There exist  $ C , M > 0 $ such that, for all $r,r_1 \in \Ra $ and all $ \eta,\eta_1 \in \Ra^{n-1} $
\begin{eqnarray}
 w (r-\log \scal{\eta}) \leq C w (r_1 - \log \scal{\eta_1})
\left( 1 + |r-r_1| + |\eta - \eta_1| \right)^M .
\label{temperateweight}
\end{eqnarray}
\end{lemm}

\noindent {\it Proof.} By Taylor's formula, $ w (x) = w (x_1) +
\int_0^1 w^{\prime} (x_1 + t (x-x_1)) \mbox{d}t (x-x_1) $ and
since
$$ w^{\prime} (x_1 + t (x-x_1)) \leq C_1  \leq
C_2 w (x_1)
$$
for all $ x,x_1 \in \Ra $ and $ t \in [0,1] $, we have $ w (x)
\leq C w (x_1) (1+|x-x_1|) $. The result then easily follows from
the fact that $ |\log \scal{\eta} - \log \scal{\eta_1}| \leq C(1 +
|\eta-\eta_1|) $. \finpreuve

\medskip

As a consequence, for all $ s \in \Ra $, $
(w(r-\log\scal{\eta}))^s $ is also a temperate weight. Hence, by
 well known pseudo-differential calculus \cite{Horm3} on $ \Ra^n $,
for all $ a \in {\mathcal S} (w_{-s}) $ and $ b \in {\mathcal
S}(w_{s}) $
\begin{eqnarray}
 a (r,y,D_r,D_y) b (r,y,D_y) = c (r,y,D_r,D_y) \label{compositionweight}
\end{eqnarray}
for some $ c \in {\mathcal S}(w_0)$ (depending continuously on $a$
and $b$). In particular, by the  Calder\`on-Vaillancourt theorem, $
c(r,y,D_r,D_y) $ is a bounded operator on $L^2$. More generally,
if $ a $ and $ b$ describe respectively bounded subsets of $
{\mathcal S}(w_{-s}) $ and $ {\mathcal S}(w_s) $, then $ c
(r,y,D_r,D_y) $ describe a bounded subset of the space of bounded
operators on $ L^2 $ (the norm of $ c (r,y,D_r,D_y) $ depends on
finitely many semi-norms of $ c $ in  ${\mathcal S}(w_0)$).
Similarly, if $a \in {\mathcal S}(w_{-s}) $, then
\begin{eqnarray} a (r,y,D_r,D_y)^* = a^{\#}(r,y,D_r,D_y) \label{adjointSws}
\end{eqnarray}
for some $ a^{\#} \in {\mathcal S}(w_{-s}) $ depending continuously on $ a $.

For $s \geq  0 $,  we introduce $ W_{s} $ as the  inverse
(unbounded if $ s \ne  0$) of $ W_{-s} $, i.e. $ W_{s} \equiv 1 $ on $
L^2 ({\mathcal K}) $ and is defined  on $ L^2
({\mathcal M} \setminus {\mathcal K}) $ as the pullback of the operator $
\widetilde{W}_{s} $ defined on $ L^2 (I \times Y) $ by
$$ (\widetilde{W}_{s} \varphi) (r,\omega) = \sum_{k \geq 0} w^{s}
(r -  \log \sqrt{ \scal{\mu_k} } ) \varphi_k (r) \psi_k (\omega) .
$$
It is clearly well defined on the dense subspace  of functions with
fast decay with respect to $r $. Then, Theorem $ \refe{theo2} $ will
clearly follow from the fact that $ W_{s} \kappa O \!p(a)
\tilde{\kappa} $ and $ \kappa O \!p(a) \tilde{\kappa} W_{s} $, defined
on $ C_c^{\infty}({\mathcal M}) $, have bounded closures on $ L^2
({\mathcal M}) $. We only consider $ W_{-s} \kappa O \!p(a)
\tilde{\kappa} $, the other case follows by adjunction, using $
\reff{adjointSws} $.

We will use a complex interpolation argument and thus
 we will need to consider $ w_{s+i\sigma} (r,\eta):= (w
 (r - \log \scal{\eta}))^{s+i\sigma} $ for $ s,\sigma \in \Ra $ (note that $ w_{s+i\sigma}
 \in {\mathcal S}(w_s) $). Since any $ a \in
{\mathcal S} (w_{-s}) $ can be written $ w_{-s} \tilde{a} $ for
some $ \tilde{a} \in {\mathcal S}(w_0) $, it is clearly enough to
show that, for all $ b \in {\mathcal S}(w_0) $, there exists $ C
> 0 $ and $ N \geq 0 $ such that
\begin{eqnarray}
 \left| \left| W_{1} \kappa O \! p (w_{-1 + i \sigma}b)
\tilde{\kappa} \varphi
 \right| \right| \leq C (1+|\sigma|)^N || \varphi ||, \qquad \forall \ \varphi \in C_c^{\infty}({\mathcal M}), \ \forall \
 \sigma \in \Ra ,  \label{principal}
\end{eqnarray}
 and that, for all $ \varphi \in C_c^{\infty} ({\mathcal M})$, there
 exists $ C_{\varphi} $ such that
\begin{eqnarray}
   \left| \left| W_{s} \kappa O \! p (w_{-s + i \sigma}b) \tilde{\kappa} \varphi
 \right| \right| \leq C_{\varphi} (1+|\sigma|)^N, \qquad  \forall \
 \sigma \in \Ra , \ \forall \ s \in [0,1]   .
\end{eqnarray}
Observing that $ W_s \scal{r}^{-1} $ is bounded, this last estimate clearly follows from the fact
that one can write
$$  W_{s} \kappa O \! p (w_{-s + i \sigma}b) \tilde{\kappa} =  W_{s}
\scal{r}^{-1} \left( \scal{r} \kappa O \! p (w_{-s + i \sigma}b)
\tilde{\kappa} \scal{r}^{-1} \right) \scal{r}
$$
and the fact that $ ||  \scal{r} \kappa O \! p (w_{-s + i \sigma}b)
\tilde{\kappa} \scal{r}^{-1} || \leq C (1 + |\sigma|)^N $, by the
Calder\`on-Vaillancourt theorem.

 We thus have to focus on $ \reff{principal} $ which we shall prove
 by using a pseudo-differential
approximation of $ W_{1} $. To that end we observe that, if $ \xi
$ is defined as in the beginning of Section $ \refe{opconjugue} $,
then
$$ w ( r - \log \scal{\mu}^{1/2} ) =  (r - \log
\scal{\mu}^{1/2} ) \xi (r - \log \scal{\mu}^{1/2})  + c (r,\mu)
$$ with $ c \in L^{\infty}( \Ra_r  \times \Ra_{\mu} )$. Thus, by choosing
 $ \chi = \chi (r) $  supported in $(r_0+2,\infty)$ such that 
$ \chi = 1 $ near 
infinity, it is easy to check that, with the notations used in
Proposition $ \refe{approxpseudo} $,
$$ W_{1} = (r - \log \scal{\Delta_{h}}^{1/2}) \chi (r) \xi (r - 
\log \scal{\Delta_{h}}^{1/2})  + B $$
for some bounded operator $ B $. Since $ || \kappa O \! p (w_{-s + i
\sigma}b) \tilde{\kappa}  || \leq C (1 + |\sigma|)^N $,  the contribution of 
$B$ to $ \reff{principal} $ is clear.
It remains to prove the following
\begin{prop} For all $ b \in {\mathcal S}(w_0) $, there exist $C > 0$ and $ N
  > 0 $ such that, fopr all $ \sigma \in \Ra $,
$$ \left| \left| (r - \log \scal{\Delta_{h}}^{1/2}) \chi (r) \xi (r - 
\log \scal{\Delta_{h}}^{1/2}) \kappa O \! p (w_{- 1 + i \sigma}b) 
\tilde{\kappa} \right| \right| \leq C (1+|\sigma|)^N . $$
\end{prop}

\noindent {\it Proof.} Observe first that $   (r - 
\log \scal{\Delta_{h}}^{1/2}) \kappa O \! p (w_{- 1 + i \sigma}b) 
\tilde{\kappa}  $ reads
\begin{eqnarray}
 \kappa O \! p (  (r-\log \scal{p_h}^{1/2})  w_{- 1 + i \sigma}b) 
\tilde{\kappa} + B_{\sigma}   \label{quatre}
\end{eqnarray}
for some bounded operator $ B_{\sigma}  $ with norm bounded by $ C
(1+|\sigma|)^N $. This follows from  the Calder\`on-Vaillancourt theorem and
the pseudo-differential expansion of $\log \scal{\Delta_h}^{1/2}  $ given
by Proposition $  \refe{approxpseudo}  $. We next insert the partition of unit
$$ 1 = \xi (r-\log \scal{p_h}^{1/2}) + (1-\xi)  (r-\log \scal{p_h}^{1/2})    $$
in front of the symbol of the first term of $ \reff{quatre}  $. Since 
 $ (r-\log \scal{p_h}^{1/2}  )  \times  \xi (r-\log \scal{p_h}^{1/2}) \times
 w_{-1+ i  \sigma}  $ belongs to $ {\mathcal S}(w_0)  $, the contribution of this
 term is clear. Thus we  are left with the study of
\begin{eqnarray}
 \chi (r) \xi (r-\log \scal{\Delta_h}^{1/2} )  
\kappa O \! p \left(  (r-\log \scal{p_h}^{1/2}) (1-\xi) (r-\log
  \scal{p_h}^{1/2} )  w_{- 1 + i \sigma}b \right) 
\tilde{\kappa} . \label{quatre2}    
\end{eqnarray}
We observe that
$$  \kappa \times  (r-\log \scal{p_h}^{1/2}) (1-\xi) (r-\log
  \scal{p_h}^{1/2} )  w_{- 1 + i \sigma} \in S^{\epsilon}  $$
for all  $ \epsilon > 0  $, since $ r \leq \log  \scal{p_h}^{1/2} + C  $ on the
support of this symbol. Then, by using the pseudo-differential expansion of
 $ \xi (r - \log \scal{\Delta_h}^{1/2}  )  $, we see that $ \reff{quatre2}  $
 reads
$$ \chi (r)  
\kappa O \! p \left(  (r-\log \scal{p_h}^{1/2}) (1-\xi) (r-\log
  \scal{p_h}^{1/2} )  \xi (r-\log
  \scal{p_h}^{1/2} )  w_{- 1 + i \sigma}b \right) 
\tilde{\kappa} + \widetilde{B}_{\sigma}  $$
with $ \widetilde{B}_{\sigma} $ similar to $ B_{\sigma} $. The symbol of the
first term belongs to $ {\mathcal S} (w_0)  $ since  $  r - \log \scal{p_h}^{1/2}  $ must be bounded on its
support and the  Calder\`on-Vaillancourt theorem completes the proof. \finpreuve

\appendix

\section{Operators on the real line} \label{dimension1}

\setcounter{equation}{0} If we consider a function $ a \in
C^{\infty} (\Ra,\Ra) $, with $ a^{\prime} $  bounded, then the
flow $ \gamma_t $, i.e. the solution to
\begin{eqnarray}
 \dot{\gamma}_t = a (\gamma_t), \qquad \gamma_0 (r) = r, \label{equadiff}
\end{eqnarray}
  is
well defined on $ \Ra_t \times \Ra_r $. For each $ t $, $ \gamma_t
$ is a $ C^{\infty} $ diffeomorphism on $ \Ra $ and it is easy to
check that
\begin{eqnarray}
 U_t \varphi : = (\partial_r \gamma_t)^{1/2} \varphi \circ
\gamma_t \label{formegroupe}
\end{eqnarray}
 defines a strongly continuous unitary group $ (U_t)_{t
\in \Ra} $ on $ L^2 (\Ra)$  whose generator, i.e. the operator $A
$ such that $U_t = e^{itA} $ for all $t$, is a selfadjoint
realization of the differential operator
$$  \frac{a(r)D_r + D_r a(r)}{2} = a (r)D_r + \frac{a^{\prime}(r)}{2i} , $$
meaning that, restricted to $ C_0^{\infty}(\Ra) $, $A $ acts as
the operator above. Indeed, according to Stone's Theorem
\cite{RS1}, the domain of $ A $, $D(A)$, is the set of $ \varphi
\in L^2 (\Ra) $ such that $ U_t \varphi $ is strongly
differentiable at $ t = 0 $, thus it clearly contains $
C_0^{\infty}(\Ra) $, which is moreover invariant by $ U_t $. This
also easily implies that $ A  $ acts on elements of its domain in
the distributions sense.

 It is worth noticing as well that, if, for some $ R $, $ a (r) = 0 $ for $ r
\leq R $, then $ \gamma_t (r) = r  $ for $ r \leq R $ and thus $
U_t $ acts as the identity on $ L^2 (-\infty,R) $. Moreover, if $
\zeta \in C_0^{1} (\Ra) $ and $ \varphi \in D (A) $, then $ \zeta
\varphi \in D (A) $ since $ U_t (\zeta \varphi) = \zeta \circ
\gamma_t U_t \varphi $ is easily seen to be strongly
differentiable at $ t = 0 $ and we have
\begin{eqnarray}
 A (\zeta \varphi) = \zeta A \varphi - i
a \zeta^{\prime}  \varphi . \label{cutoffloin}
\end{eqnarray}
Of course, it is not hard to deduce from this property that the
subspace of $ D (A) $ consisting of compactly supported elements
is dense in $ D (A) $ for the graph norm.

We want to show that $ C_0^{\infty} (\Ra) $ is also a core for $ A
$ and thus consider $ \theta_{\epsilon} (r) = \epsilon^{-1} \theta
(r / \epsilon) $ with $ \theta \in C_0^{\infty} (-1,1) $ such that
$ \int_{\Ra} \theta = 1 $. A simple calculation shows that
$$ U_t (\varphi \ast \theta_{\epsilon}) = K_{t,\epsilon} U_t \varphi  $$
where  $ K_{t,\epsilon} $ is the operator with kernel
$$ \kappa_{t,\epsilon} (r,r^{\prime}) = \left( \partial_r \gamma_t (r) \right)^{1/2}
\left( \partial_r \gamma_t (r^{\prime}) \right)^{1/2}
\theta_{\epsilon}(\gamma_t(r) - \gamma_t (r^{\prime})) . $$ Note
that this operator is bounded on $ L^2 (\Ra) $ in view of the
following well known Schur's Lemma which we recall since we will
use it extensively.
\begin{lemm}[Schur] If $ j (x,y) $ is a measurable function on $
\Ra^{2d} $ such that
$$  \emph{ess-}\sup_{y \in \Ra} \int |j (x,y)| \ \emph{d}x \leq C , \qquad \emph{ess-} \sup_{x \in \Ra}
\int |j(x,y)| \ \emph{d}y \leq C $$ then the operator $J$ with
kernel $ j $ is bounded on $ L^2 (\Ra^d) $ and $ ||J|| \leq C $.
\end{lemm}
Since $ K_{0,\epsilon} \varphi = \varphi \ast \theta_{\epsilon} $,
we have
\begin{eqnarray}
 \frac{U_t (\varphi \ast \theta_{\epsilon}) - \varphi \ast
\theta_{\epsilon} }{it} -
 \left( \frac{U_t \varphi - \varphi}{it} \right) \ast \theta_{\epsilon} =
 \frac{K_{t,\epsilon} - K_{0,\epsilon}}{it} (U_t \varphi) .
 \label{approxdiff}
\end{eqnarray}
In order to estimate the right hand side, we start with a few
remarks. Note first that we have
\begin{eqnarray}
 ||\partial_r \gamma_t||_{\infty} \leq
e^{||a^{\prime}||_{\infty}|t|} , \qquad || \partial_r \partial_t
\gamma_t ||_{\infty} \leq ||a^{\prime}||_{\infty}
e^{||a^{\prime}||_{\infty}|t|} . \label{Gronwall}
\end{eqnarray}
 The first estimate is obtained
by applying $ \partial_r $ to $ \reff{equadiff} $ and using
Gronwall's lemma. The second one then follows from the first one.
This implies in particular the existence of some $ t_0 $,
depending only on $ ||a^{\prime}||_{\infty} $, such that $ ||
\partial_r \gamma_t - 1 ||_{\infty} \leq 1/2 $ for $ |t| \leq t_0
$. Differentiating $ \reff{equadiff} $ twice with respect to $ r$
yields $ \partial_t^2
\partial_r \gamma_t = a(\gamma_t)a^{\prime \prime}(\gamma_t) \partial_r \gamma_t + a^{\prime}(\gamma_t)^2
\partial_r \gamma_t $ and thus, if $ a
a^{\prime \prime} $ is bounded,
\begin{eqnarray}
 || \partial_t^2 \partial_r
\gamma_t ||_{\infty} \leq \frac{3}{2} \left( ||a a^{\prime
\prime}||_{\infty} +||a^{\prime}||^2_{\infty} \right) , \qquad |t|
\leq t_0 . \label{Gronwall2}
\end{eqnarray}
 Thus,  if $ J_{\epsilon} $ denotes the operator with kernel $
\partial_t \kappa_{t,\epsilon}|_{t=0} (r,r^{\prime})
 $ that is
\begin{eqnarray}
 \frac{1}{2} \left( a^{\prime} (r) + a^{\prime} (r^{\prime})
\right) \theta_{\epsilon} (r-r^{\prime}) + \left( a (r) - a
(r^{\prime}) \right) \theta_{\epsilon}^{\prime}(r-r^{\prime}) ,
\label{kappa}
\end{eqnarray}
 then Taylor's formula combined with Schur's Lemma show that
\begin{eqnarray}
 \left| \left| K_{t,\epsilon} - K_{0,\epsilon} - t J_{\epsilon}
\right| \right| \leq C_{\epsilon}t^2, \qquad |t| \leq t_0
\label{estimeepsilon}
\end{eqnarray}
 for some $ C_{\epsilon} $ depending only on $ \theta_{\epsilon} $, $
||a^{\prime}||_{\infty} $ and $ ||a a^{\prime \prime}||_{\infty} $
(recall that $ t_0 $ depends only on $ ||a^{\prime}||_{\infty} $
as well). Since $ J_{\epsilon} $ is a bounded operator (with norm
uniformly bounded by $ ||a^{\prime}||_{\infty} \int |r
\theta^{\prime}(r)| + |\theta(r)| \mbox{d}r $),  $
\reff{approxdiff} $ and $ \reff{estimeepsilon} $ show that if $
\varphi \in D (A) $ then $ \varphi \ast \theta_{\epsilon} \in D
(A) $ and $ A (\varphi \ast \theta_{\epsilon}) = (A \varphi) \ast
\theta_{\epsilon} - i J_{\epsilon} \varphi $. Furthermore $
J_{\epsilon} \rightarrow 0 $ strongly as $ \epsilon \rightarrow 0
$ for it is uniformly bounded and  $ J_{\epsilon} \psi \rightarrow
0 $ for all $ \psi \in C_0^{\infty} (\Ra) $. All this shows that,
for any $ \varphi \in D (A) $,
\begin{eqnarray}
|| A (\varphi \ast \theta_{\epsilon}) - (A \varphi) \ast
\theta_{\epsilon} || \leq C || \varphi ||, \qquad  A (\varphi \ast
\theta_{\epsilon}) - (A \varphi) \ast \theta_{\epsilon}
\rightarrow 0 , \qquad \epsilon \rightarrow 0 ,
\label{convolution}
\end{eqnarray}
with $ C $ independent of $ \epsilon $, depending only on $ ||
a^{\prime}||_{\infty} $. In particular, $ \reff{cutoffloin} $ and
$ \reff{convolution} $ imply easily that $ C_0^{\infty}(\Ra) $ is
a core for $ A $.

\section{Proof of Proposition $ {\bf \refe{approxpseudo}} $} \label{functional}
\setcounter{equation}{0} We start with some reductions. We may
clearly write $ g (r,\mu) $ as $ g_1(r,\mu) (i + \mu) $ with $ g_1
\in S^{-1} $ hence by studying $ g_1(r,\Delta_{h}) $ instead of $ g
(r,\Delta_{h}) $ we can assume that $g \in S^m $ with $ m < 0 $.
Note that the composition by $ \Delta_{h} + i $ on the right of $
\reff{compright} $ doesn't cause any trouble in view of $
\reff{regularisation} $, $ \reff{comDr} $ and of the standard
composition rules for pseudo-differential operators. Furthermore,
by positivity of $ \Delta_{h} $, we have $ g (r, \Delta_{h}) = g_2 (r
, \Delta_{h} + 1) $ for some $ g_2 \in S^m $ which we can assume to
be supported in $ [ 1/2 , \infty ) $. This support property will
be useful to consider Mellin transforms below.

By the standard procedure for the calculus of a parametrix of the
resolvent of an elliptic operator on a closed manifold
\cite{Seel1}, there are symbols $
q_{-2}(y,\eta,z),q_{-3}(y,\eta,z),\cdots $ of the form
\begin{eqnarray}
q_{-2} = (p_h -z)^{-1}, \qquad q_{-2-j} = \sum_{1 \leq l \leq 2 j}
d_{jl} (p_h - z)^{-l-1}, \qquad j \geq 1
\end{eqnarray}
  such that, for all $ N $ large enough,
$$ \theta (\Delta_{h} - z)^{-1}  - \Psi^* \left( (\Psi_* \theta)
\sum_{j = 0}^{N} q_{-2-j}(y,D_y,z)   \right) \Psi_* \
\tilde{\theta} =  M_N (z).
$$
Here $ M_N (z) $ is bounded from $ H^{\kappa}  $ to $ H^{\kappa+N}
 $ for all $ \kappa $, $ H^{\kappa} =  H^{\kappa} (Y)$ being the standard
Sobolev space on $Y$
 and $ d_{jl} $ are polynomials in $ \eta $  of degree $ 2 j - l $,
which are independent of $ z $ and linear combinations of products
of derivatives of the full symbol of $ \Delta_{h} $ in the chart we
consider. Furthermore, for all $ \kappa $ and $ N $,
there exist $ C $ and $ \gamma $ such that
$$ || M_N (z) ||_{H^{\kappa} \rightarrow H^{\kappa+N}} \leq
C  \frac{\scal{z}^{\gamma}}{|\mbox{Im}z|^{\gamma+1}} .
$$
We now repeat the arguments of \cite{HeRo1}. For each $s$ such
that $ \mbox{Re} s < 0 $, we choose a contour $ \Gamma_s $
surrounding $ [1/2 , + \infty ) $ on which $ \scal{z} /
|\mbox{Im}z| $ is bounded, and by Cauchy formula we get
$$ \theta (\Delta_{h} + 1)^{s}  - \Psi^* \left( (\Psi_* \theta)
\sum_{j = 0}^{N} a_{j}(y,D_y,s)   \right) \Psi_* \ \tilde{\theta}
= \frac{i}{(2 \pi)} \int_{\Gamma_s} z^s  M_N (z) \ \mbox{d}z $$
with $ a_j (s) = \sum_{1 \leq l \leq 2 j} (-1)^l d_{jl}  s (s-1)
\cdots (s-l+1) (p_h + 1)^{s-l} / l!  $ if $ j \geq 1 $ and $ a_0
(s) = (p_h + 1)^{s} $. As in \cite{HeRo1}, we choose the contour
so that, if $ \mbox{Re} \ s < 0  $ is fixed,
$$  \left| \left|  \int_{\Gamma_s}   z^s  M_N (z) \ \mbox{d}z \right| \right|_{H^{\kappa} \rightarrow
H^{\kappa+N}} \leq C_{\mathrm{Re} s, \kappa , N} \scal{\mbox{Im}
s}^{\gamma} .
$$
We then consider the Mellin transform $ \mathrm{M}[g_2](r,s) : =
\int_{0}^{\infty} \mu^{s-1} g_2 (r,\mu) \ \mbox{d}\mu $.   Note
that it is well defined for $ s < - m $ (recall that $ m < 0 $),
since $ g_2 $ is supported in $ [1/2,\infty ) $, and that it
decays fast at infinity with respect to $ |\mbox{Im}s| $, for
fixed $ \mbox{Re} \ s $. It is then easy to check that
$$  \int_{ \mathrm{Re} s - i \infty}^{\mathrm{Re} s + i \infty}  \mathrm{M}[g_2] (r,s)
\left(  \int_{\Gamma_s} z^s  M_N (z) \ \mbox{d}z \right) \
\mbox{d}s \in C^{\infty}(\Ra_r , {\mathcal
L}(H^{\kappa},H^{\kappa+N})) ,
$$
so, by Mellin's inversion formula, i.e. $ g_2 (r, \mu) = (2 i \pi)^{-1}
\int_{\mathrm{Re} s = \mathrm{const}} \mathrm{M}[g_2](r,s)
\mu^{-s} \ \mbox{d}s $, and by setting
$$ {\mathcal R}^{\theta,Y}_N (r) = \theta g(r,\Delta_{h} )  - \Psi^* \left( (\Psi_* \theta)
g(p_h ) + (\Psi_* \theta) \sum_{j = 1}^{N} \sum_{l} (-1)^l d_{jl}
\partial_{\mu}^{l}g(r,p_h) / l! \right) \Psi_* \ \tilde{\theta} ,$$ we
get
$$ \sup_{r > r_0 }  \left| \left| D_r^k {\mathcal R}_N^{\theta, Y } (r)
\right| \right|_{H^{\kappa} \rightarrow H^{\kappa+N} } < \infty, \qquad  \forall \ k . $$
The latter easily follows from the boundedness of the derivatives
of $ g $ (or $g_2$) with respect to $ r$. In order to prove $
\reff{regularisation} $, with $ N $ replaced by $ N/8 $ (which can
be assumed to be an integer), we first remark that $ {\mathcal
R}^{\theta}_N $ is defined on generators of $ L^2 (I ) \otimes L^2 (Y) $
by
$$  {\mathcal R}^{\theta}_N  (\varphi_k \otimes \psi_k)(r,\omega)  =
\varphi_k(r)\left( {\mathcal R}^{\theta,Y}_{N} (r) \psi_k \right)
(\omega) $$ with $ \varphi_k \in L^2 (I) $. We then note that, by
writing $ \psi_k = (\mu_k + i)^{-N/4} (\Delta_{h} + i )^{N/4} \psi_k
$, we have, for $ j , l \leq N/8 $,
$$ \left| \left| \Delta_{h}^j {\mathcal R}^{\theta}_N \Delta_y^l (\varphi_k \otimes \psi_k)
\right| \right|_{L^2(I \times Y)} \leq C_N \scal{\mu_k}^{-N/4}
||\varphi_k ||_{L^2 (I)} \sup_{r > r_0} || {\mathcal R}_{
N}^{\theta,Y}(r) ||_{H^{- 3N/4} \rightarrow H^{N/4}}
$$
and thus, if $ N $ is large enough so that  $ \sum_k
\scal{\mu_k}^{-N/2} < \infty $, Parseval's formula yields
$$   \left| \left| \Delta_{h}^j {\mathcal R}^{\theta}_N \Delta_{h}^l ( \sum_k \varphi_k \otimes \psi_k)
\right| \right|_{L^2(I \times Y)} \leq C_N  \left( \sum_k
||\varphi_k ||^2_{L^2 (I)} \right)^{1/2} . $$ This proves $
\reff{regularisation} $. The proof of $ \reff{comDr} $ is similar.
\finpreuve


\begin{thebibliography}{99}



\bibitem{ABG1} {\sc W. Amrein, A. Boutet de Monvel, V. Georgescu}, {\it
    $C_0$-Groups, Commutators methods and spectral theory of N-body
    hamiltonians},
Birkh\"auser (1996).

\bibitem{BrPe1} {\sc V. Bruneau, V. Petkov},
{\it Semiclassical resolvent estimates for trapping
perturbations}, Commun. Math. Phys. 213, no. 2, 413-432 (2000).

\bibitem{Bouc4} {\sc J.M. Bouclet}, {\it Spectral distributions for long range perturbations},
J. Funct. Analysis, 212, no. 2, 431-471 (2004). 

\bibitem{Bouc5}  {\sc \name}, {\it Generalized scattering
phases  for asymptotically hyperbolic manifolds}, CRAS, 338 no. 9, 685-688 (2004).

\bibitem{Bouc8} {\sc \name}, {\it A Weyl law for asymptotically hyperbolic
    manifolds}, in preparation.


\bibitem{CaVo1} {\sc F. Cardoso, G. Vodev},
{\it Uniform estimates of the resolvent of the Laplace-Beltrami
operator on infinite volume Riemannian manifolds. II}, Ann. Henri
Poincar\'e 3, 673-691 (2002).

\bibitem{ChPi} {\sc J. Chazarain, A. Piriou},
{\it Introduction \`a la th\'eorie des \'equations aux d\'eriv\'ees partielles
lin\'eaires},
Gauthier-Villars, Paris, (1981).

\bibitem{DiSj1} {\sc M. Dimassi, J. Sj\"ostrand},
{\it Spectral asymptotics in the semi-classical limit}, London
Mathematical Society Lecture Note Series, 268. Cambridge
University Press, (1999).

\bibitem{FrHi1} {\sc R. G. Froese, P. D. Hislop},
{\it Spectral analysis of second-order elliptic operators on
noncompact manifolds}, Duke Math. J. 58, no. 1, 103-129 (1989).


\bibitem{GeGe1}  {\sc V. Georgescu, C. G\'erard}, {\it On the virial theorem
    in quantum mechanics}, Commun. Math. Phys. 208, no. 2, 275-281 (1999).


\bibitem{HeRo1} {\sc B. Helffer, D. Robert},
{\it Calcul fonctionnel par la transformation de Mellin et
op\'erateurs admissibles}, J. Funct. Analysis 53, 246-268 (1983).

\bibitem{Horm3}  {\sc L.  H\"ormander}, {\it  The analysis  of  linear partial
differential operators III}, Springer-Verlag (1985).

\bibitem{IsKi1} {\sc H. Isozaki, H. Kitada}, {\it Modified wave operators with
  time independent modifiers}, J. Fac. Sci., University of Tokyo, 
Section I A 32, 77-104 (1985).

\bibitem{IsKi2} {\sc \name}, {\it Microlocal resolvent estimates for 
2-body Schr\"odinger operators},
J. Funct. Analysis  57, no. 3, 270-300 (1984), and {\it Erratum}
J. Funct. Analysis  62, no. 2, 336 (1985).





\bibitem{Melr0} {\sc R.B. Melrose}, {\it Geometric scattering
theory}, Stanford lectures, Cambridge Univ. Press (1995).


\bibitem{Mour1} {\sc E. Mourre}, {\it Absence of singular continuous spectrum for certain
self-adjoint operators}, Commun. Math. Phys. 78, 391-408 (1981).

\bibitem{Mour2} {\sc \name}, {\it Op\'erateurs conjugu\'es et propri\'et\'es
 de propagation}, Commun. Math. Phys. 91, no. 2, 279-300 (1983).

\bibitem{PSS1} {\sc P. Perry, I. M. Sigal, B. Simon}, {\it Spectral analysis
of  N-body Schr\"odinger operators}, Ann. Math. 114, no. 3,
519-567 (1981).

\bibitem{RS1} {\sc M. Reed, B. Simon}, {\it Modern methods in mathematical
physics, vol. I}, Academic Press.

\bibitem{Robe1} {\sc D. Robert}, {\it Relative time delay for perturbations of
  elliptic operators and semiclassical asymptotics}, J. Funct. Analysis 126,
  no. 1, 36-82 (1994).

\bibitem{Robe2} {\sc \name}, {\it On the Weyl formula for obstacles}, 
 Partial differential equations and mathematical physics,  
 Progr. Nonlinear Differential Equations Appl., 21, Birkhäuser, 264-285 (1996). 


\bibitem{Seel1} {\sc R.T. Seeley},
{\it  Complex powers of an elliptic operator}, Singular integrals
(Proc. Sympos. Pure Math. Chicago, III 1966) A.M.S. R.I., 288-307
(1967).


\bibitem{Vode2} {\sc G. Vodev}, {\it Local energy decay of solutions
to the wave equation for non trapping metrics}, Ark. Math. 42,
379-397 (2004).

\end{thebibliography}
\end{document}